\documentclass{amsart}

\usepackage[ left=2.8cm, right=2.8cm]{geometry}

\usepackage{setspace}
\usepackage{xcolor}

\theoremstyle{plain}
\newtheorem{thm}{Theorem}[section] 
\newtheorem{rem}[thm]{Remark}

\usepackage{amssymb}
\usepackage{amsmath}
\usepackage{lineno}
\usepackage{lipsum}
\usepackage{amsfonts}
\usepackage{graphicx, color}
\usepackage{subfigure}
\usepackage{float}    
\usepackage{comment}
\usepackage{epstopdf}
\usepackage{dsfont}
\usepackage{todonotes}
\usepackage{ulem}

\newcommand{\REV}[1]{{\color{black}{#1}}}

\begin{document}

\title[Flocking dynamics with topological interactions]{\Large Numerical modeling of flocking dynamics with topological interactions}

\author[M. Menci]{Marta Menci}\address[M.M.]{Università Campus Bio-Medico di Roma, Rome, Italy}
\email{m.menci@unicampus.it}

\author[T. Paul]{Thierry Paul}
\address[T.P.]{CNRS  Laboratoire Ypatia des Sciences Mathématiques (LYSM),  Rome, Italy }
\email{thierry.paul@cnrs.fr}

\author[S. Rossi]{Stefano Rossi}
\address[S. R.]{ETH, Zurich, Switzerland }
\email{stefano.rossi@math.ethz.chr}

\author[T. Tenna]{Tommaso Tenna}\address[T. T.]{Laboratoire J.A.Dieudonné, Université Côte d’Azur, Nice, France \& Università degli Studi di Roma La Sapienza, Rome, Italy}
\email{tommaso.tenna@univ-cotedazur.fr}

\begin{abstract}
In this paper, we propose a numerical investigation of topological interactions in flocking dynamics. 
Starting from a microscopic description of the phenomena, mesoscopic and macroscopic models have been previously derived under specific assumptions. We explore the role of topological interactions by describing the convergence speed to consensus in both microscopic and macroscopic dynamics, considering different forms of topological interactions. Additionally, we compare mesoscopic and macroscopic dynamics for monokinetic and non-monokinetic initial data. Finally, we illustrate with some simulations in one- and two-dimensional domains the sensitive dependence of solutions on initial conditions, including the case where the system exhibits two solutions starting with the same initial data.

\end{abstract}

\maketitle

\section{Introduction}
\label{sec1}
Collective dynamics have been widely investigated over the last years in different fields of research, due to the great variety of living and non-living systems exhibiting such complex behaviors. From flocks of birds, schools of fish to bacterial colonies and robotic swarms, coordinated movements of individuals interacting within a group give raise to the so called \textit{emergent behaviors}. Understanding the underlying principles of self-organized behaviors has emerged as a fundamental challenge in applied mathematics, aiming at providing the framework to describe, predict, and analyze the basic mechanisms determining the observed collective dynamics.
A crucial aspect of mathematical modeling is selecting the appropriate scale for the description for the phenomena under investigation. In traditional microscopic, agent-based approaches, the dynamics of each individual in the group are tracked.
When the number of agents increases, mesoscopic and macroscopic models can be considered as efficient alternatives to pure microscopic ones. This is mainly due to the fact that microscopic models, focusing on the dynamics of each individual, become computationally demanding for large groups.
Conversely, mesoscopic and macroscopic approaches allow to describe collective behaviors using aggregated variables such as probabilistic distribution of agents, macroscopic density and average velocity field (see \cite{TE_2024} for an updated reference paper on the derivation from microscopic to macroscopic models). 
 Several mathematical models have been proposed in the literature to model \textit{flocking} behaviors, see \cite{VZ} for an introduction to the field. Typically displayed by flocks of birds, flocking is used to define the coordinated motion of agents reaching the same velocity when moving towards the same direction. The majority of the models, starting from the seminal papers \cite{cucker2007emergent, cucker2007mathematics}, assumes interactions involving metric quantities (\textit{metric interactions}), where individuals adjust their movement in response to the other members of the group, and the interaction depends on the relative Euclidean distance.
However, increasing empirical evidence suggests that in many biological systems, interactions are \textit{topological}, meaning that individuals interact with a fixed number of neighbors, regardless of their absolute distances \cite{Ballerini,Bialek}. As a consequence, several mathematical models adopting this novel paradigm have appeared, ranging from microscopic to macroscopic ones, see e.g. \cite{Haskovec2013, BD2016, ST_topo, KD, CMPB21, AF2024}.
This physically non-metric nature of interactions plays a crucial role for the stability, robustness and efficiency of the observed motion, affecting not only the local organization of individuals but also the emergent macroscopic patterns.
\subsection{\textit{Topological interactions models: from microscopic to macroscopic scale}}
In this paper we focus on the topological model of the Cucker-Smale type firstly proposed in \cite{Haskovec2013}. Starting from the microscopic scale, the model reads 
\begin{equation}
\label{sec2:CSmodel_pij}
    \begin{cases}
        \dot X_i(t) =V_i(t) \\
        \dot V_i(t)= \displaystyle \frac{1}{N} \sum_{j=1}^N p_{i,j}(V_j(t) - V_i(t)),
    \end{cases}
    \quad i\in\{1,\dots,N\}.
\end{equation}
Here $X_i\in \mathbb{R}^d$, $V_i\in \mathbb{R}^d$ denote positions and velocities of $N$ indistinguishable interacting agents in a general $d\ge 1$ dimensional space. A key factor of the model is the so-called communication weights $p_{i,j}$, modeling the kind of considered interactions. 

 In the original version of the Cucker-Smale  model \cite{cucker2007emergent, cucker2007mathematics}, the interaction $p_{i,j}$ 
between two agents $i,j \in \{1, \dots, N\}$ 
depends on the metric distance between them, 
i.e. $p_{i,j}= p_{i,j}(|X_i-X_j|)$.

In the topological framework, $p_{i, j}$ are chosen as
\begin{equation}
\label{sec2:Kfunction}
    p_{i,j}:= K\left( M^N_{i,j}\right),
\end{equation}
where $K:[0,1] \to \mathbb{R}^+$ is a non increasing function. The quantity $M^N_{i,j}$ is defined as
\begin{equation}
\label{toprank}
    M^N_{i,j}:= \frac{1}{N}\sum_{h=1}^N \chi_{\bar B(X_i,|X_i-X_j|)}(X_h),
\end{equation}
where $\chi_{\bar B(X_i,|X_i-X_j|)}$ denotes the characteristic function of the closed ball centered in $X_i$ with radius $|X_i-X_j|$.
The quantity $M^N_{i,j}$ hence counts the number of agents inside the ball $\bar B(X_i, |X_i-X_j|)$, and it is known as \textit{topological rank} between two agents $i,j \in \{1,\dots, N\}$.
Note that non-decreasing property of function $K$ reflects the fact that the farther apart two agents are, also in a topological sense, the less influence one has on the other. In the next section, specific choices of non-increasing function $K$ are performed, see the following \eqref{function_K1},\eqref{function_K2}.

Starting from the microscopic model \eqref{sec2:CSmodel_pij},  
the kinetic and hydrodynamic equations
have been first formally derived in \cite{Haskovec2013}.
At kinetic level, the Vlasov-type equation related to \eqref{sec2:CSmodel_pij} reads as
\begin{equation}
\begin{aligned}
    \label{sec2:Vlasov}
        \partial_t f_t(x,v)+ &v \cdot \nabla_x f_t(x,v)\\
        &\quad = \nabla_v \cdot 
        \left( f_t(x,v)
        \int K\left(M[\nu^0_t](x,|x-y|)\right)(v-w) df_t(y,w) \right),
    \end{aligned}
\end{equation}
where $f_t$ denotes a distribution function of agents at time $t$ depending on $(x,v) \in \mathbb{R}^d \times \mathbb{R}^d$. We have introduced the following notation for the zero-th and first order moments of $f_t= f_t(x,v) $, respectively
\begin{equation}\label{nu0}
\nu^0_t(x) :=\int f_t(x,v) \, dv,
\end{equation}
\begin{equation}\label{nu1}
\nu^1_t(x) :=\int v f_t(x,v) \, dv
\end{equation}
for any $t\ge0$ and
\begin{equation}
\label{sec2:M}
    M[\nu^0_t](x,|x-y|):= \int_{\bar B(x, |x-y|)} d\nu^0_t
\end{equation}
the mass inside the ball $\bar B(x,|x-y|)$.

At the macroscopic scale, the following pressureless Euler-type system is formally obtained in \cite{Haskovec2013}:
\begin{equation}
  \label{eulereq}
  \left\{
  \begin{aligned}
    &\partial_t \rho_t(x)+ \nabla_x \cdot (\rho_t(x) u_t(x))=0, \\
    &\partial_t u_t(x) + (u_t(x) \cdot \nabla_x)u_t(x) =
        \int_{\mathbb{R}^d}
     K\Big(M[\rho_t](x,|x-y|)\Big)(u_t(y)-u_t(x))\rho_t(y) d y, \\
  \end{aligned}
  \right.
\end{equation}
under a monokinetic assumption on the kinetic distribution function, namely
\begin{equation}
    f_t(x,v) = \rho_t(x) \, \delta ( v - u_t(x)).
\end{equation}
Here 
 $(\rho_t(x), u_t(x)) :
[0,T]\times \mathbb{R}^d \to \mathbb{R} \times \mathbb{R}^d$, $T>0$, denote density and velocity field, respectively. 

\begin{rem}
    The quantities $\nu^0_t$ and $\nu^1_t$ denote mass and momentum of the distribution function $f_t$ at the kinetic level. Under a monokinetic assumption, the kinetic moments and the macroscopic moments coincide as 
    \begin{equation}
        \nu_t^0(x) = \rho_t (x), \qquad \nu_t^1(x) = \rho_t(x)\,u_t(x).
    \end{equation}
\end{rem}

We remark that the derivation in \cite{Haskovec2013} is only formal. Indeed the rigorous proof of the mean-field limit and propagation of chaos results for the topological kinetic equation \eqref{sec2:Vlasov}, under the assumption of Lipschitz continuity for the interaction function $K$ in \eqref{sec2:Kfunction}, has been presented in \cite{BCR21, BPR24}. The results on propagation of chaos have been further extended to the case of general $K:[0,1] \to [0,1]$ non-increasing functions in \cite{MPR25}.
From the analytical point of view, local and global well-posedeness results for the multi-dimensional case of the Euler system can be found in \cite{RS, LRS}, and rigorous results have been also obtained for stochastic models of topological type \cite{DP19, DPR23}.

From a numerical perspective, the literature is still quite poor and several issues need to be tackled.
In this work we present a numerical investigation of the topological models in \eqref{sec2:CSmodel_pij}-\eqref{eulereq} at different scales. In particular, we focus on different choices of function $K$ in \eqref{sec2:Kfunction}, which rules topological interactions, studying how this affects the convergence to consensus. Moreover, we investigate the similarities and differences between the mesoscopic and macroscopic dynamics both within and outside the monokinetic regime. Finally, we illustrate with examples in one- and two-dimensional scenarios the sensitive dependence of the solutions on the initial data.

The paper is organized as follows: in Section \ref{sec:Numerical_approx} we detail the numerical schemes implemented to approximate the solutions at the kinetic and macroscopic level. The main results of our study are reported in  Section \ref{sec:Numerical_tests}, where different scenarios, at different scales, in 1D and 2D setting are discussed, including the case with non uniqueness of solutions. Conclusions and future perspective on the topic conclude the paper in Section \ref{sec:Conclusion}.

\section{Numerical Approximation}\label{sec:Numerical_approx}
In this section we describe the numerical schemes implemented to approximate the solutions of the models at different scales.
The microscopic model is numerically approximated using classical implicit Euler scheme both in the one and two dimensional case, avoiding restrictions on the time step due to explicit integration. The kinetic and the macroscopic models are approximated using finite volume schemes, as detailed in the following.
The time discretization on $[0,T]$ is performed introducing a discretization time step $\Delta t$, which is chosen according to a CFL stability condition for the kinetic and the macroscopic equation.
\subsection{Numerical scheme for the kinetic model}\label{sec:numscheme_kinetic}
We detail the construction of the scheme in the one dimensional case, using a classical upwind scheme both in the spatial and in the velocity direction. Let us consider uniform spatial cells $\mathcal{C}_i = \left(x_{i-1/2}, x_{i+1/2}\right)$, where $x_i$ are equally distributed points and $x_{i+1/2}=x_i+\Delta x / 2$. Analogously, we define uniform velocity cells  $\mathcal{V}_k = \left(v_{k-1/2}, v_{k+1/2}\right)$, where $v_k$ are equally distributed points and $v_{k+1/2}=v_k+\Delta v / 2$. Let us consider the kinetic equation \eqref{sec2:Vlasov} in the one-dimensional case rewritten as
\begin{equation}
    \partial_t f + \partial_x F(f) + \partial_v G(f) = 0,
\end{equation}
where 
\begin{equation}
\label{functions_kinetic_numerical}
    F (f)  = v\,f, \qquad G (f) = - f \int K \left(M [\nu^0_t] (x, |x-y|) \right)\,(v-w) df(y,w) .
\end{equation}
The fully discrete upwind scheme reads as
\begin{equation}
\label{numerical_kinetic}
    f^{n+1}_{i,k} = f^n_{i,k} - \frac{\Delta t}{\Delta x} \left( F^n_{i+1/2,k} -F^n_{i-1/2,k} \right) - \frac{\Delta t}{\Delta v} \left( G^n_{i,k+1/2} -G^n_{i,k-1/2} \right),
\end{equation}
where $F^n_{i+1/2,k}$ and $G^n_{i,k+1/2}$ are first-order upwind numerical fluxes. 
In particular, we can explicitly define
\begin{equation}
    F^n_{i+1/2, k} = v_k\,f^n_{i+1/2,k} = \begin{cases}
        v_k\,f^n_{i,k}, \quad &\text{if } v_k > 0,\\
       v_k\,f^n_{i+1,k}, \quad &\text{if } v_k \leq 0,
    \end{cases}
\end{equation}
and by introducing $\xi_{i,k+1/2}$ as a quadrature approximation of the integral term in \eqref{functions_kinetic_numerical} at $(x_i,v_{k+1/2})$, we can also define 
\begin{equation}
    G^n_{i,k+1/2} = \xi_{i,k+1/2}\,f^n_{i,k+1/2} =       
    \begin{cases}
        \xi_{i,k+1/2}\,f^n_{i,k}, \quad &\text{if } \xi_{i,k+1/2} > 0,\\
        \xi_{i,k+1/2}\,f^n_{i,k+1}, \quad &\text{if } \xi_{i,k+1/2} \leq 0.
    \end{cases}
\end{equation}
The choice of a first order scheme relies in its positivity-preserving property, which is fundamental for a probability distribution function. We briefly describe the derivation of the CFL condition to guarantee the positivity-preserving of the scheme, inspired by \cite{Tan2017}.
\subsubsection{Positivity preserving property}
It is possible to rewrite \eqref{numerical_kinetic} as
\begin{multline}
\label{Positivity_Preserving_Eq}
    f^{n+1}_{i,k} = \frac{1}{2} \left(f^n_{i,k} - \frac{2\Delta t}{\Delta x} v_k\,f^n_{i+1/2,k} - \frac{2\Delta t}{\Delta v} \xi_{i,k+1/2}\,f^n_{i,k+1/2}\right) +\\ \frac{1}{2} \left(f^n_{i,k} + \frac{2\Delta t}{\Delta x} v_k\,f^n_{i-1/2,k} + \frac{2\Delta t}{\Delta v} \xi_{i,k-1/2}\,f^n_{i,k-1/2}\right).
\end{multline}
Then, if $\xi_{i,k+1/2} \leq 0$, we have
\begin{multline}
    f^n_{i,k} - \frac{2\Delta t}{\Delta x} v_k\,f^n_{i+1/2,k} - \frac{2\Delta t}{\Delta v} \xi_{i,k+1/2}\,f^n_{i,k+1/2} \\ = f^n_{i,k} - \frac{2\Delta t}{\Delta x} v_k\,f^n_{i+1,k} - \frac{2\Delta t}{\Delta v} \xi_{i,k+1/2}\,f^n_{i,k+1} > 0, \qquad \text{ if} \, v_k \leq 0,
\end{multline}
and
\begin{multline}
    f^n_{i,k} - \frac{2\Delta t}{\Delta x} v_k\,f^n_{i+1/2,k} - \frac{2\Delta t}{\Delta v} \xi_{i,k+1/2}\,f^n_{i,k+1/2} \\ = \left(1- \frac{2\Delta t}{\Delta x} v_k\right)\,f^n_{i,k} - \frac{2\Delta t}{\Delta v} \xi_{i,k+1/2}\,f^n_{i,k+1} > 0, \qquad \text{ if} \, v_k>0,
\end{multline}
under the CFL condition $\Delta t < \displaystyle \frac{\Delta x}{2\,\max_k|v_k|}$. In the same spirit, if $\xi_{i,k+1/2} > 0$
\begin{multline}
    f^n_{i,k} - \frac{2\Delta t}{\Delta x} v_k\,f^n_{i+1/2,k} - \frac{2\Delta t}{\Delta v} \xi_{i,k+1/2}\,f^n_{i,k+1/2} \\ = \left( 1 - \frac{2\Delta t}{\Delta v} \xi_{i,k+1/2} \right) \, f^n_{i,k} - \frac{2\Delta t}{\Delta x} v_k\,f^n_{i,k+1}  > 0, \qquad \text{ if} \, v_k \leq 0,
\end{multline}
under the CFL condition $\Delta t < \displaystyle \frac{\Delta v}{2\,\max_{i,k} |\xi_{i,k+1/2}|}$. Again we have
\begin{multline}
    f^n_{i,k} - \frac{2\Delta t}{\Delta x} v_k\,f^n_{i+1/2,k} - \frac{2\Delta t}{\Delta v} \xi_{i,k+1/2}\,f^n_{i,k+1/2} \\ = \left( 1 - \frac{2\Delta t}{\Delta x} v_k - \frac{2\Delta t}{\Delta v} \xi_{i,k+1/2} \right) \, f^n_{i,k} > 0, \qquad \text{ if} \, v_k>0,
\end{multline}
under the CFL condition $\Delta t < \displaystyle \frac{\Delta v\,\Delta x}{2\,\max_{i,k} \left(|\Delta v \, v_k + \Delta x \, \xi_{i,k+1/2}|\right)}$.\\[12pt]
In the same way, the second term in \eqref{Positivity_Preserving_Eq} is positive under the same condition. We conclude that the first order scheme under the CFL condition
\begin{equation}
\label{CFL_kinetic}
    \Delta t < \displaystyle \frac{\Delta v\,\Delta x}{2\,\max_{i,k} \left(|\Delta v \, v_k + \Delta x \, \xi_{i,k+1/2}|\right)},
\end{equation}
is positivity preserving, namely  $f^{n+1}_{i,k} > 0$ if $f^n_{i,k} > 0$ for all $i$ and $k$.

\subsection{Numerical scheme for the macroscopic model}\label{sec:schema_macro}

Let us describe the numerical scheme for the approximation of the macroscopic system \eqref{eulereq} in the one- and two-dimensional cases. We present the details of the approximation in the two-dimensional case, since the one-dimensional case can be easily derived.
\REV{In the two-dimensional case denoting $u = (u_1, u_2)$, $\nabla=(\partial_x, \partial_y)$ and omitting the subscript $t$ for simplicity of notation, system \eqref{eulereq} reads as
\begin{equation}
  \label{eulereq_2D_system}
  \left\{
  \begin{aligned}
    &\partial_t \rho+ \partial_x (\rho u_1) + \partial_y (\rho u_2) =0, \\
    &\partial_t u_1 + u_1\partial_x u_1 + u_2\partial_y u_1 = \displaystyle \int_{\mathbb{R}^2} K \Big(M[\rho](\mathbf{z},|\mathbf{z}-\mathbf{z}'|)\Big)(u_1(\mathbf{z}')-u_1(\mathbf{z}))\rho(\mathbf{z}') d \mathbf{z}',\\
    &\partial_t u_2 + u_1\partial_x u_2 + u_2\partial_y u_2 = \displaystyle \int_{\mathbb{R}^2} K \Big(M[\rho](\mathbf{z},|\mathbf{z}-\mathbf{z}'|)\Big)(u_2(\mathbf{z}')-u_2(\mathbf{z}))\rho(\mathbf{z}') d \mathbf{z}'
  \end{aligned}
  \right.
\end{equation}
where $\mathbf{z}=(x,y) \in \mathbb{R}^2$.
The numerical scheme for the approximation of the solution combines the finite volume second-order central scheme proposed by Kurganov and Tadmor \cite{KT}. A variant of this scheme has been implemented in \cite{ACH2018} for the numerical approximation of a pressureless Euler alignment system without topological interactions.
We here extend the approach proposed in \cite{KT} with the global flux technique \cite{Gascon2001, Ciallella2023, Kurganov2023} to treat non-conservative terms.
In this spirit, we rewrite system \eqref{eulereq_2D_system} in a quasi-conservative form as
\begin{equation}
\label{def_w_S}
    \partial_t w + \partial_x \tilde{A}(w) + \partial_y \tilde{B}(w) = S(w),
\end{equation}
where 
\begin{equation}
    w = \begin{pmatrix}
        \rho\\
         u_1\\
         u_2
    \end{pmatrix},
    \qquad 
S(w) = \begin{pmatrix}
        0\\[10pt]
        \displaystyle \int_{\mathbb{R}^2} K \Big(M[\rho](\mathbf{z},|\mathbf{z}-\mathbf{z}'|)\Big)(u_1(\mathbf{z}')-u_1(\mathbf{z}))\rho(\mathbf{z}') d \mathbf{z}'\\[10pt]
        \displaystyle \int_{\mathbb{R}^2} K \Big(M[\rho](\mathbf{z},|\mathbf{z}-\mathbf{z}'|)\Big)(u_2(\mathbf{z}')-u_2(\mathbf{z}))\rho(\mathbf{z}') d \mathbf{z}'
    \end{pmatrix}.
\end{equation}
The fluxes are defined by
\begin{equation}
\tilde{A}(w) = \begin{pmatrix}
        \rho u_1\\
        \frac{1}{2} u_1^2\\
        \displaystyle u_1 u_2 + \int_{\hat{x}}^x N_1(\xi,y)\, d\xi
    \end{pmatrix},
    \qquad
\tilde{B} (w) = \begin{pmatrix}
        \rho u_2\\
        \displaystyle u_1 u_2 + \int_{\hat{y}}^y N_2(x,\eta)\, d\eta\\
         \frac{1}{2} u_2^2
    \end{pmatrix}
\end{equation}
with 
\begin{equation}
    N_1(x,y) = -u_2(x,y)\partial_x u_1(x,y), \qquad N_2 (x,y) = - u_1(x,y)\partial_y u_2(x,y)
\end{equation}
and $\hat{x}$, $\hat{y}$ arbitrary values.
}
Let us consider uniform discretization in space, with cells $\REV{\mathcal{C}}_{i,j}= \left(x_{i-1/2}, x_{i+1/2}\right) \times \left(y_{j-1/2}, y_{j+1/2}\right)$ of volume $\Delta x \times \Delta y$, whose centers are given by $(x_i, y_j)$. The semi-discrete version of our scheme in the two-dimensional case reads as
\begin{equation}
\label{scheme_macro_semidiscrete}
    \frac{d}{dt} \bar{w}_{i,j} = -\frac{1}{\Delta x} \left( \REV{\mathcal{F}}^x_{i+\frac{1}{2}, j} (t) - \REV{\mathcal{F}}^x_{i-\frac{1}{2}, j}(t) \right)  - \frac{1}{\Delta y} \left( \REV{\mathcal{F}}^y_{i, j+\frac{1}{2}} (t) - \REV{\mathcal{F}}^y_{i, j-\frac{1}{2}} (t) \right) + S_{i,j}(t),
\end{equation}
where $\bar{w}_{ij}$ represents the cell average of $w$ in $\REV{\mathcal{C}}_{i,j}$ \REV{and $\mathcal{F}^x$ (or $\mathcal{F}^y$) is the numerical flux corresponding to $\tilde{A}$ (or $\tilde{B}$ respectively)}.
To obtain a second-order scheme we introduce a piecewise linear reconstruction for $w$ of the form
\begin{equation}
\label{polynomial_reconstruction}
    w_{i,j}(x,y) = \bar{w}_{i,j} + (w_x)_{i,j}\,(x-x_i) + (w_y)_{i,j}\,(y-y_j).
\end{equation}
The derivatives $w_x$ and $w_y$ are approximated by using a minmod limiter \cite{Sweby1984} to avoid oscillations in the reconstruction procedure, namely
\begin{equation*}
    (w_x)_{i,j} = \text{minmod} \, \left( \theta \frac{\bar{w}_{i+1,j}-\bar{w}_{i,j}}{\Delta x}, \theta \frac{\bar{w}_{i+1,j}-\bar{w}_{i-1,j}}{2\,\Delta x}, \theta \frac{\bar{w}_{i,j}-\bar{w}_{i-1,j}}{\Delta x} \right),
\end{equation*}
\begin{equation*}
    (w_y)_{i,j} = \text{minmod} \, \left( \theta \frac{\bar{w}_{i,j+1}-\bar{w}_{i,j}}{\Delta y}, \theta \frac{\bar{w}_{i,j+1}-\bar{w}_{i,j-1}}{2\,\Delta y}, \theta \frac{\bar{w}_{i,j}-\bar{w}_{i,j-1}}{\Delta y} \right),
\end{equation*}
where $\theta \in [1,2]$. We can now easily define the reconstructed values at the interfaces $w_{i,j \pm \frac{1}{2}}$ and $w_{i\pm \frac{1}{2}, j }$ through the expression \eqref{polynomial_reconstruction} and the characteristic speeds at the interfaces. We define the following quantities
\begin{equation*}
     w^N_{i,j}=w_{i,j}(x_i,y_{j+1/2}), \qquad  w^S_{i,j}:=w_{i,j}(x_i,y_{j-1/2}), \qquad  w^E_{i,j}:=w_{i,j}(x_{i+1/2},y_j), \qquad
     w^W_{i,j}:=w_{i,j}(x_{i-1/2},y_j).
\end{equation*}

In particular $a_{i+1/2,j}^\pm$ will denote the largest (lowest) characteristic speed of $\REV{\tilde{A}}$ along the interface $i+1/2$ and $b_{i, j+1/2}^\pm$ the largest (lowest) characteristic speed of $\REV{\tilde{B}}$ along the interface $j+1/2$.\\
The numerical fluxes read as
\begin{equation*}
    \REV{\mathcal{F}}^x_{i+\frac{1}{2}, j} = \frac{a^+_{i+\frac{1}{2},j}\, \tilde{A}(\REV{w^E_{i,j}}) - a^-_{i+\frac{1}{2},j}\, \tilde{A}(\REV{w^W_{i+1,j}})}{a^+_{i+\frac{1}{2},j} - a^-_{i+\frac{1}{2},j}} + \frac{a^+_{i+\frac{1}{2},j} \, a^-_{i+\frac{1}{2},j}}{a^+_{i+\frac{1}{2},j} - a^-_{i+\frac{1}{2},j}} \left(\REV{w^W_{i+1,j} - w^E_{i,j}} \right)
\end{equation*}
and
\begin{equation*}
    \REV{\mathcal{F}}^y_{i, j+\frac{1}{2}} = \frac{b^+_{i,j+\frac{1}{2}}\, \tilde{B}(\REV{w^N_{i,j}}) - b^-_{i,j+\frac{1}{2},}\, \tilde{B}(\REV{w^S_{i,j+1}})}{b^+_{i,j+\frac{1}{2}} - b^-_{i,j+\frac{1}{2}}} + \frac{b^+_{i,j+\frac{1}{2}} \, b^-_{i,j+\frac{1}{2}}}{b^+_{i,j+\frac{1}{2}} - b^-_{i,j+\frac{1}{2}}} \left(\REV{w^S_{i,j+1}- w^N_{i,j}} \right).
\end{equation*}
The integral term $S_{i,j}$ is approximated by using a nested quadrature rule (the first one for the computation of the quantity $M$ and the second one for the full integral).

Specifically, the approximation of the second component of the source term in \eqref{def_w_S}, here denoted as
$S_{i,j}^{\,(2)}$, is computed as
\begin{equation}
    S_{i,j}^{\,(2)} := \rho_{i,j} \sum_{\mathbf{z}'_{p,q} \in \mathbb{Z}} \omega_{\mathbf{z}'} K(M[\rho_\Delta](\mathbf{z}_{i,j}, ||\mathbf{z}_{i,j}-\mathbf{z}'_{p,q}||) \left( u_{p, q} - u_{i,j} \right) \rho_{p,q},
\end{equation}
where $\mathbf{z}_{i,j} := (x_i, y_j)$, $u_{i, j} \approx u_1(x_i,y_j)$ and the quantity $M$ is approximated by
\begin{equation}
    M[\rho_\Delta](\mathbf{z}_{i,j}, ||\mathbf{z}_{i,j}-\mathbf{z}'||) \approx \sum_{\mathbf{h}_{l,m} \in \mathbb{Z}} \tilde{\omega}_{\mathbf{h}} \, \rho_{l,m} \mathds{1}_{\{||\mathbf{h}_{l,m}-\mathbf{z}_{i,j}|| \leq ||\mathbf{z}_{i,j} - \mathbf{z}'||\}}.
\end{equation}
The quantities $\omega_{\mathbf{z}'}$ and $\tilde{\omega}_{\mathbf{h}}$ are the corresponding quadrature weights. 
The third component $S_{i,j}^{\,(3)}$ is computed analogously, by defining $u_{i, j} \approx u_2(x_i,y_j)$. 

The semi-discrete scheme obtained in \eqref{scheme_macro_semidiscrete} is discretized in time using classical explicit Euler method.

\section{Numerical Tests}\label{sec:Numerical_tests}
We present several numerical tests simulating microscopic, kinetic and macroscopic dynamics. 
We investigate the role of topological interactions on the overall dynamics, comparing the results obtained at different scales. Further comparison will involve kinetic and macroscopic dynamics. We conclude our analysis exploiting sensitive dependence on initial conditions in a 2D setting. 
All the tests performed have been run a laptop equipped with an Intel Core i7-1060NG7 processor and 16 GB RAM.

\subsection{Effects of Topological Interactions}

\begin{figure}[!h]
\centering
 \subfigure{\includegraphics[scale=0.55]{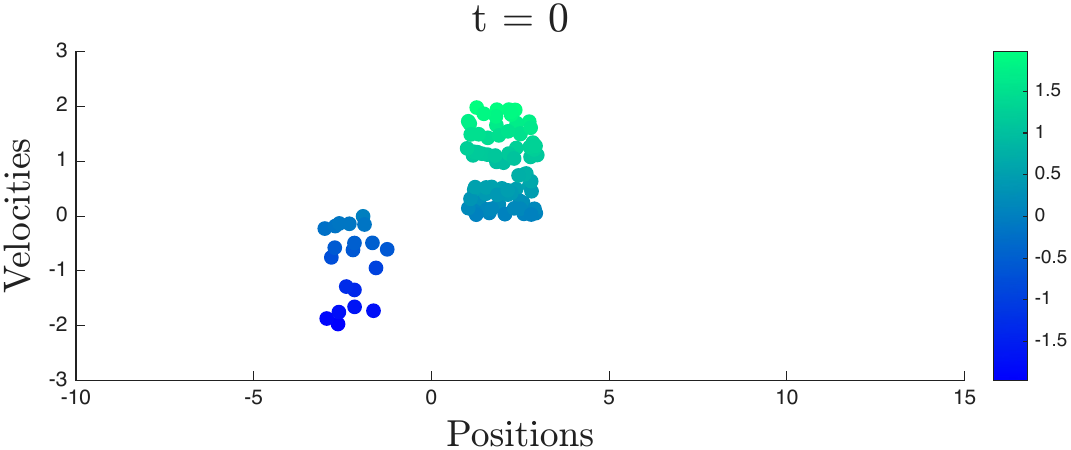}}
 \caption{ Test 1: initial condition  $N=100$ agents. Each position is marked with colors based on its initial velocity: the velocities of the agents located in $[-3,-1]$ are randomly chosen in $[-2,0]$, whereas the velocities of the agents located in $[1,3]$ are randomly chosen in $[0,2]$.}
\label{Test2gruppi_micro_CI}
\end{figure}

We investigate the role of topological interactions on the overall dynamics, comparing the results obtained for different choices of function $K$ in  \eqref{sec2:Kfunction}. 
We consider a reference microscopic scenario in 1D, with $N=100$ agents. At initial time, a smaller group of $N_1=20$ agents is randomly located in $[-3, -1]$ with initial velocity randomly assigned in $[-2, 0]$, whereas the remaining $N_2=80$ agents are located in $[1, 3]$ with initial random velocity in $[0, 2]$. Figure \ref{Test2gruppi_micro_CI} shows the initial condition plotted in the position-velocity space.  

In the following Test1-Test3 we run system \eqref{sec2:CSmodel_pij} over the time interval $[0,500]$, with time discretization step $\Delta t = 0.02$. 

\paragraph{\textit{Test 1: Study of consensus velocities for different values of $\overline{M}$}}
We here consider a non Lipschitz function, choosing $K= \mathds{1}_{[0, \overline{M}/N  ]}$. Varying the value of $\overline{M} >0$, we show the effect of variations in the number of neighbors considered for interactions. 

\begin{figure}[!h]
\centering
\subfigure[$\overline{M}=11$]{\includegraphics[scale=0.33]{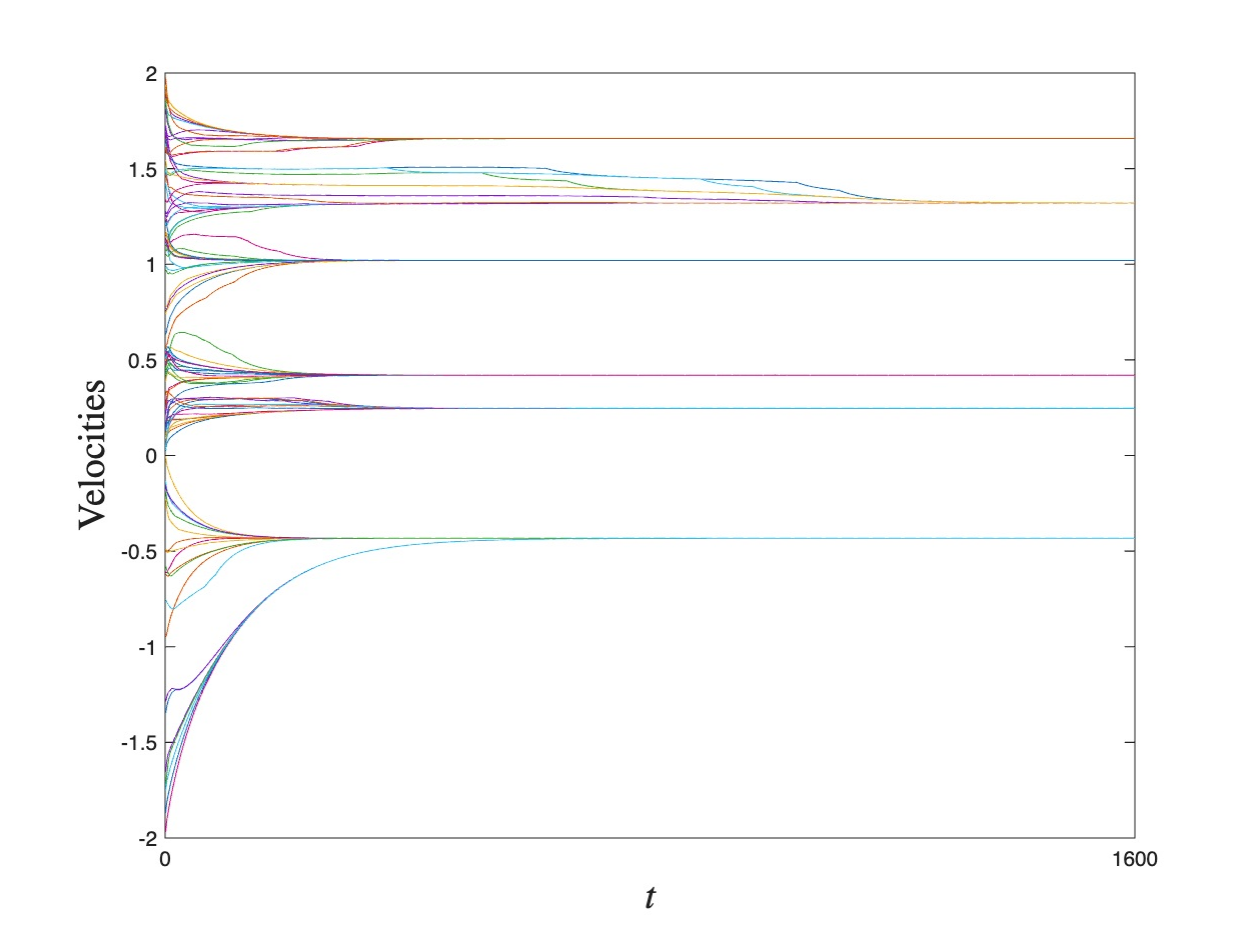}}
\subfigure[$\overline{M}=22$]{\includegraphics[scale=0.33]{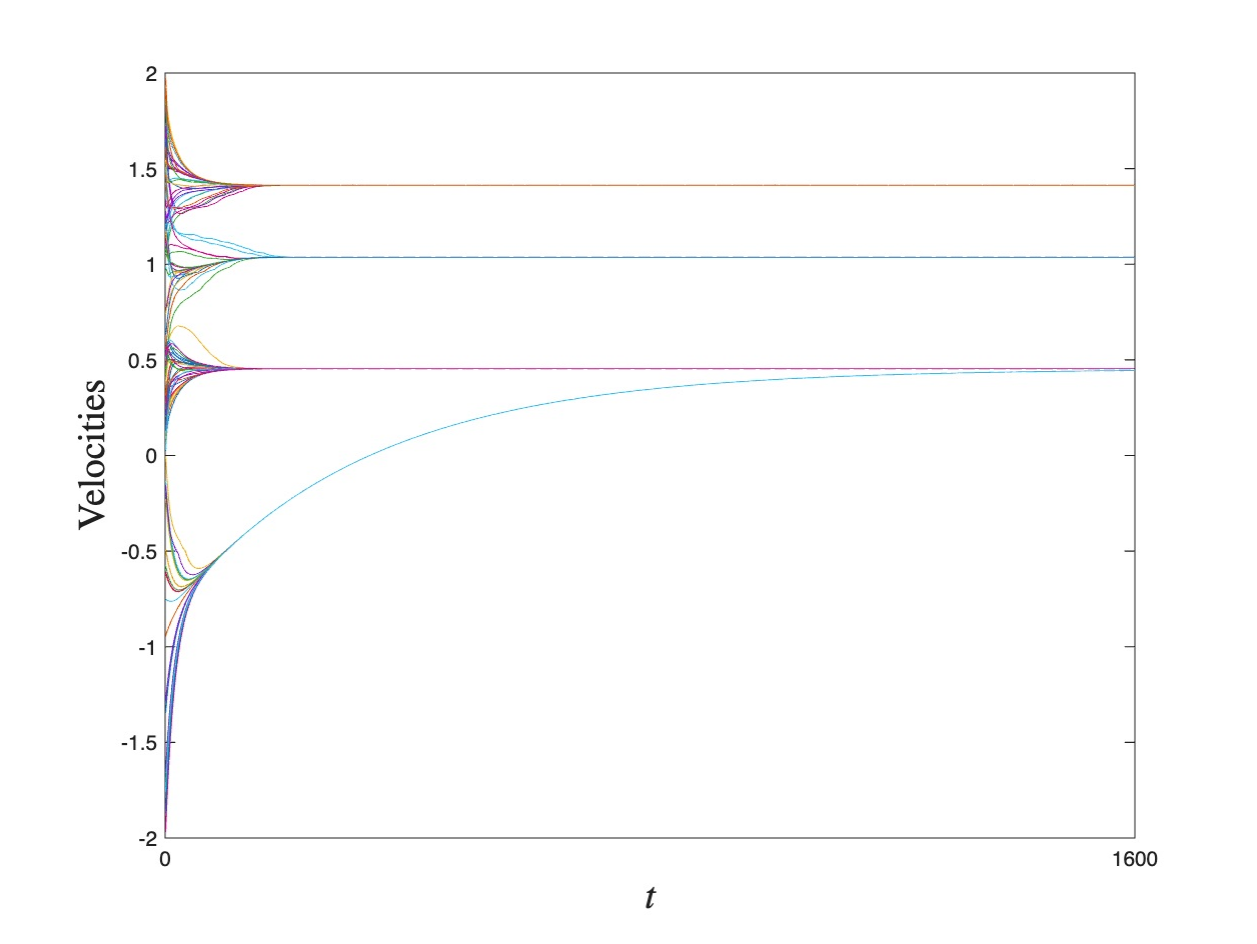}}\\
\subfigure[$\overline{M}=60$]{\includegraphics[scale=0.33]{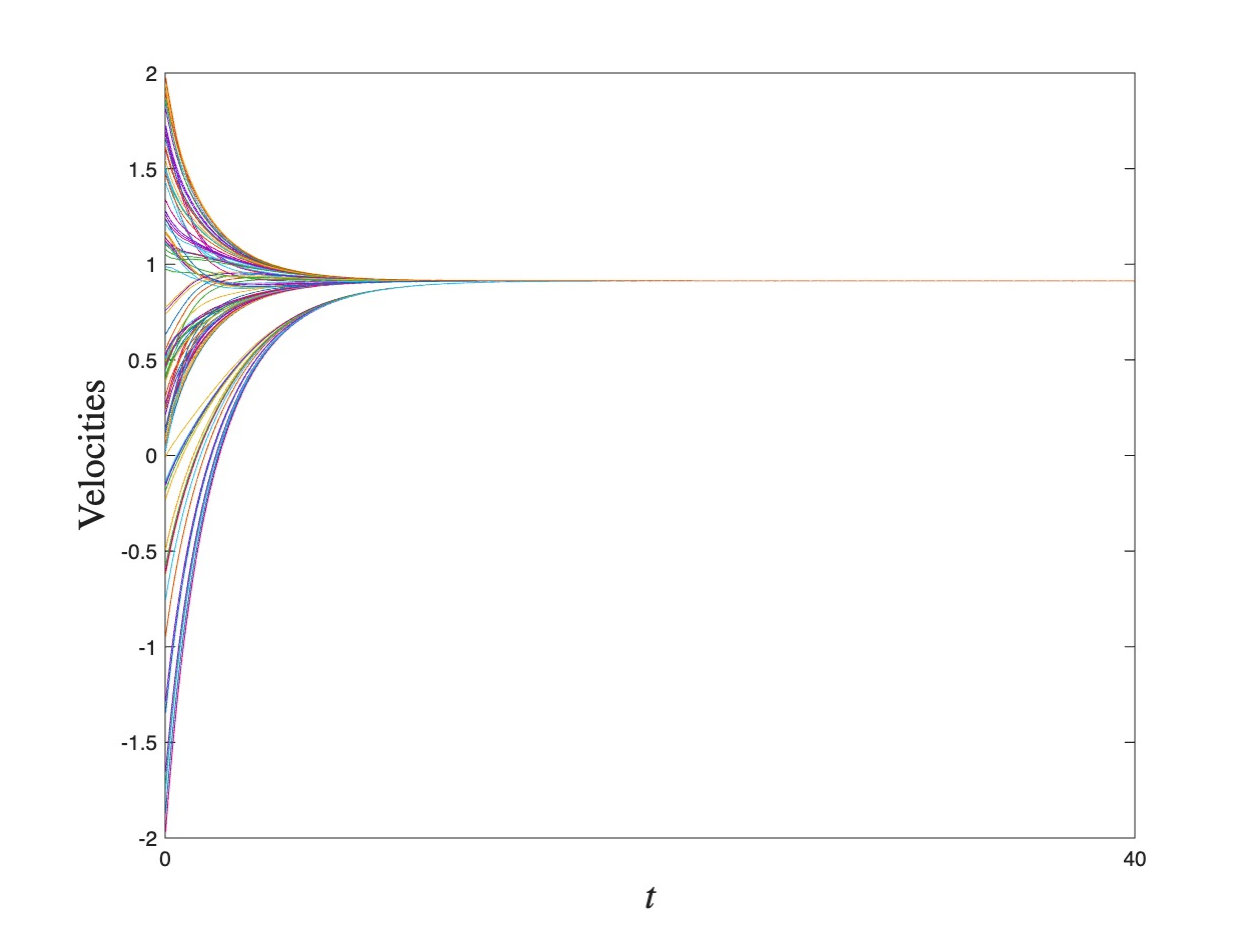}}
\subfigure[$\overline{M}=100$]{\includegraphics[scale=0.33]{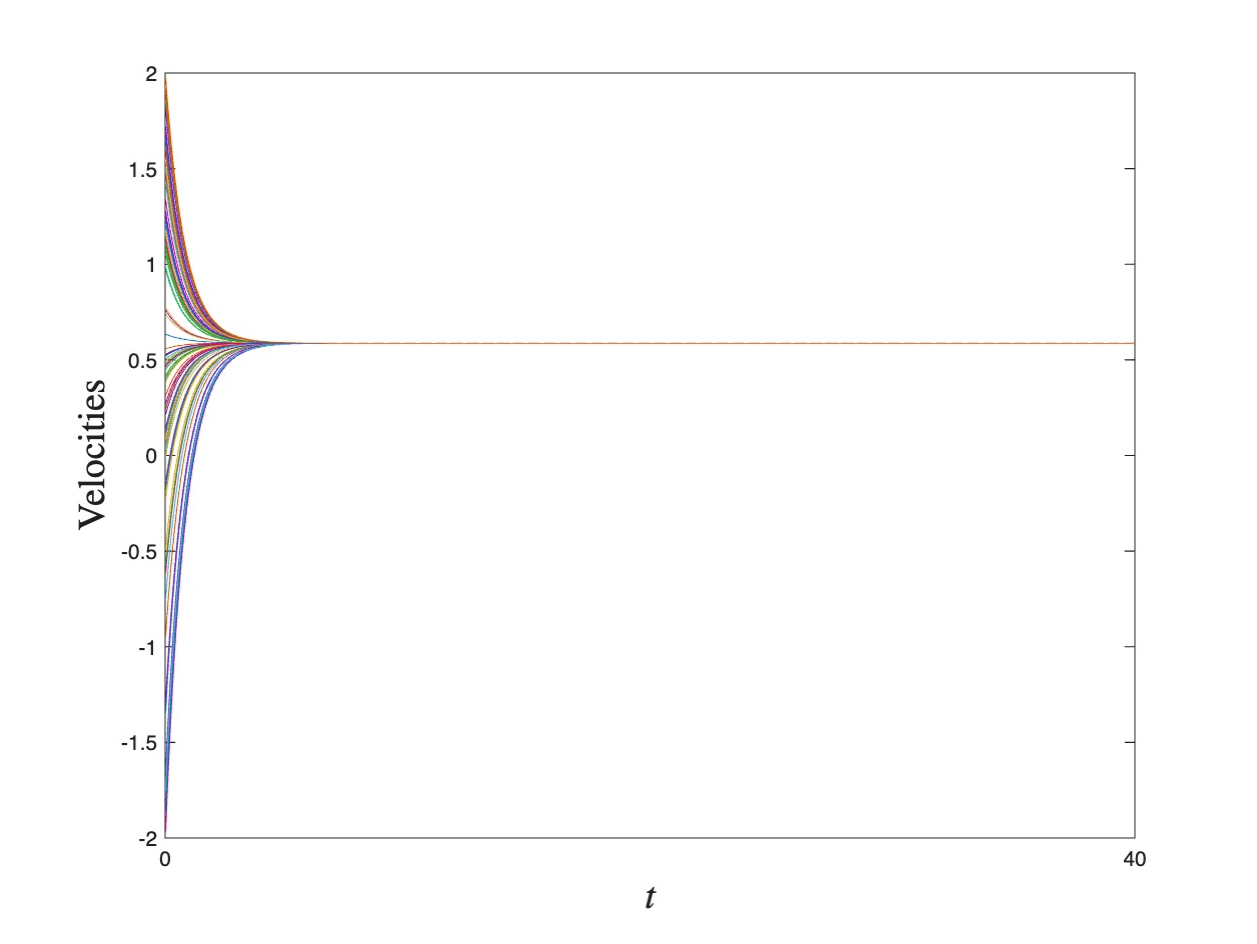}}
 \caption{ Test 1: plot of $V_i(t)$, $i=1,...,N$, solution to \eqref{sec2:CSmodel_pij} with  $K=\mathds{1}_{[0, \overline{M}/N  ]}$ for different values of $\overline{M}$.}
\label{Test_Kvicini}
\end{figure}

\begin{figure}[!h]
\centering
\subfigure[$K(x)=1-x$]{\includegraphics[scale=0.32]{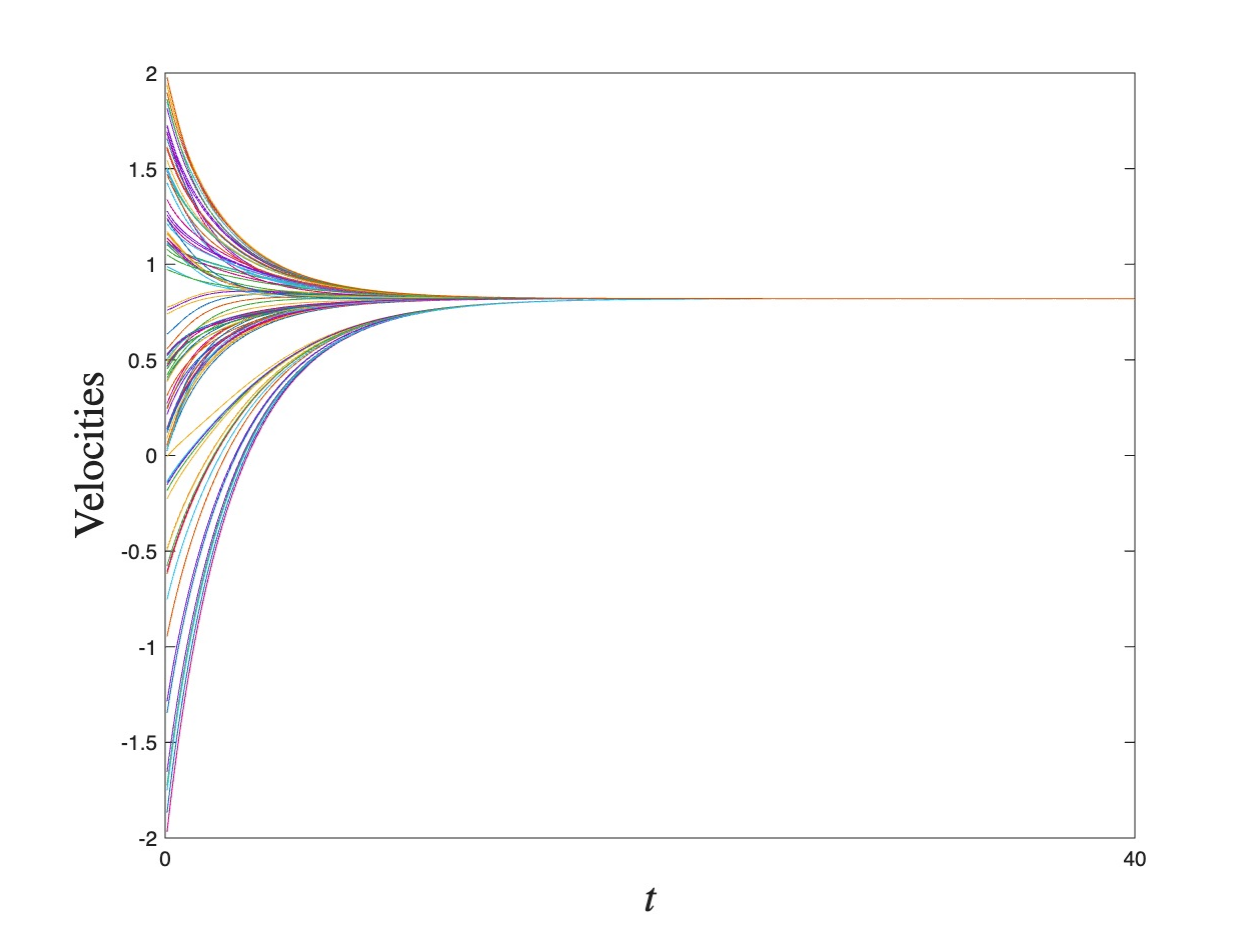}}
\subfigure[$K(x)=(1-x)^2$]{\includegraphics[scale=0.32]{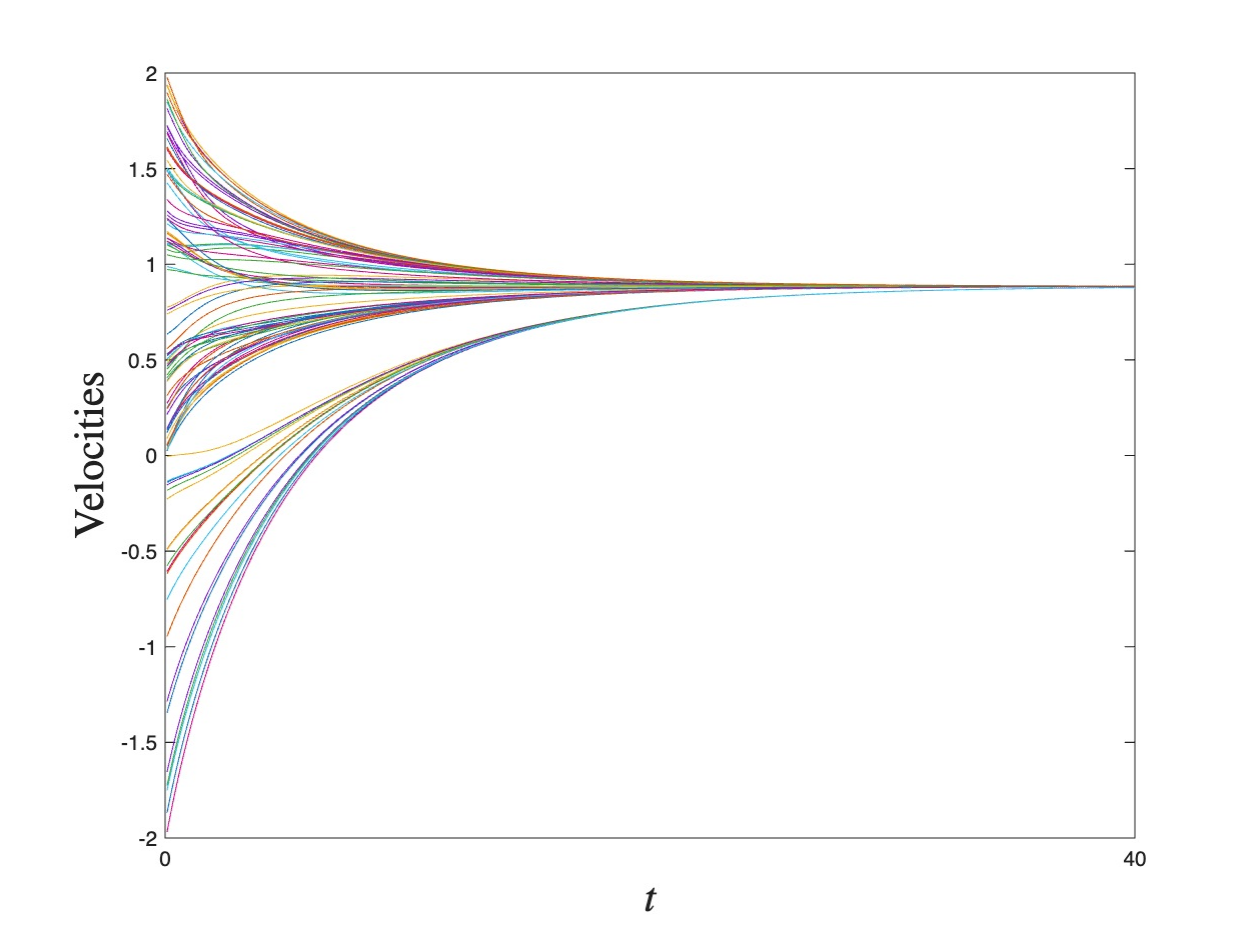}}\\
\hspace{-0.3cm}
\subfigure[Microscopic velocity]{\includegraphics[scale=0.35]{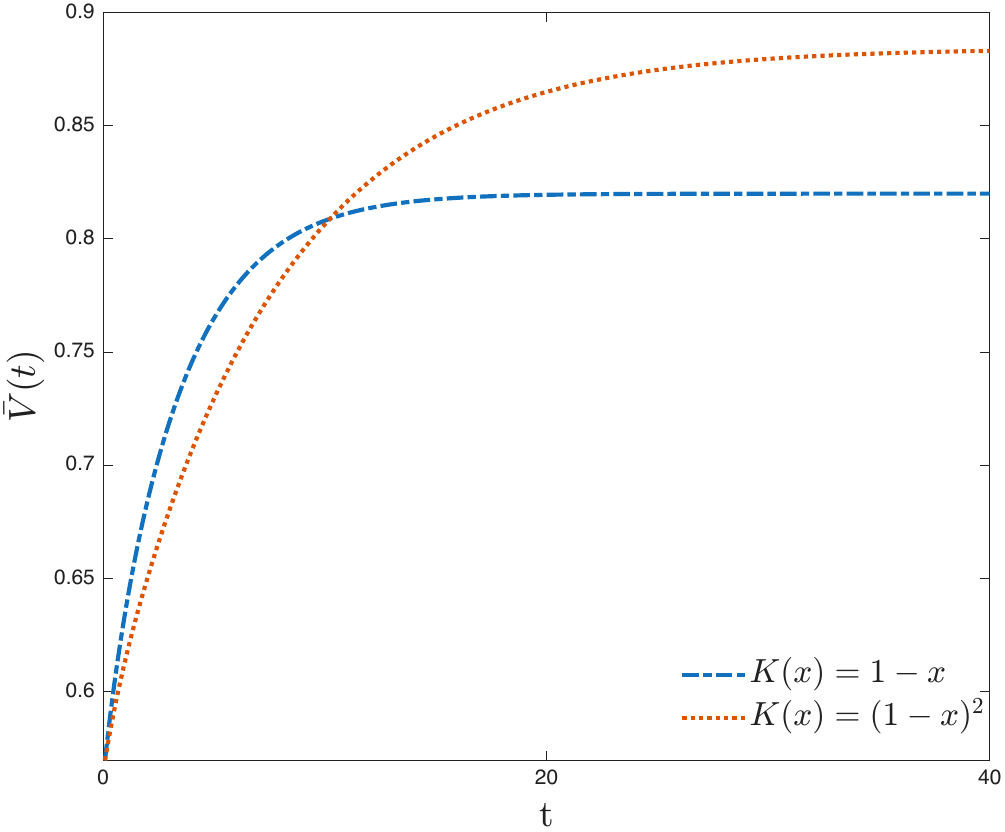}}
\hspace{0.8cm}
\subfigure[Macroscopic velocity]{\includegraphics[scale=0.35]{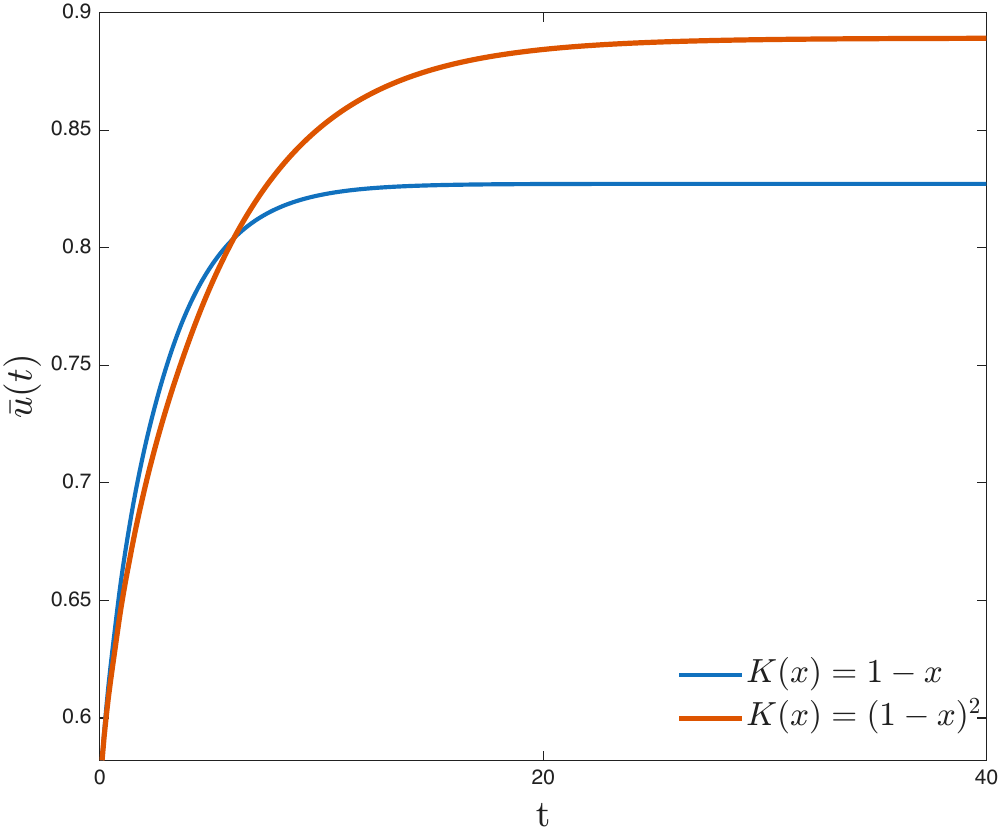}}
 \caption{ Test 2: plot of the microscopic velocity $V_i(t)$, $i=1,...,N$, solution to \eqref{sec2:CSmodel_pij} with a) $K$ linear b) $K$ convex function and the comparison between c) $\bar{V}(t)$ and d) $\bar{u}(t)$, for different $K$ functions.
 }
\label{Test_K_micromacro}
\end{figure}

Figure \ref{Test_Kvicini} shows the velocity evolution of the agents for four different values of $\overline{M}$, which are representative for the simulated scenario. Indeed, choosing $\overline{M}< N_1$ implies that agents in the smaller group do not interact with the remaining ones, and different consensus velocities are reached within the groups (Figure \ref{Test_Kvicini} (a)).  In order to let the two groups communicating, we consider $\overline{M}> N_1$ (Figure \ref{Test_Kvicini} (b)).
We observe that the number of consensus velocities decreases as $\overline{M}$ increases, since interacting with a higher number of neighbors reduces the formation of cluster dynamics.
Assuming $\overline{M}=N$ we recover an all-to-all type of interaction. Comparing the results obtained for $\overline{M}=60$, we observe that, keeping $N$ fixed, the value of $\overline{M}$ influences the time at which flocking is reached (see Figure \ref{Test_Kvicini} (c)-(d)). This is due to the fact, keeping $N$ fixed, the process of achieving consensus is faster increasing the number of neighbors each agent considers for interactions. 
\paragraph{\textit{Test 2: Study of consensus velocities for different choices of $K$}}
In Test 2 we compare the velocity of convergence to a flocking state for two different $K:[0,1] \to [0,1]$, both Lipschitz and decreasing functions:
\begin{subequations}
\begin{align}   
    \label{function_K1} &K_1 (x) = 1-x,\\[10pt]
    \label{function_K2} &K_2(x)=(1-x)^2.
\end{align}
\end{subequations}
Figure \ref{Test_K_micromacro} (a)-(b) shows the results obtained for $K_1$ and $K_2$, respectively. Function $K_2$ decays to zero faster with respect to $K_1$, reflecting in a weaker interactions among agents. This also implies a slower convergence to the flocking state.
The obtained results can be better visualized in Figure \ref{Test_K_micromacro} (c), showing the evolution in time of the microscopic mean velocity $\bar{V}$, for the two different choices of $K$. 

Moreover, we compare the results obtained at the microscopic scale with the macroscopic one.
In order to link microscopic initial data to macroscopic ones, different approaches can be adopted, see \cite{Silverman,Terrell} for a complete introduction to seminal methods in the field. Our choice allows us to introduce the initial data for the kinetic scale, that will be used in Test 5, creating a full path, from microscopic to macroscopic scale.

Starting from the initial data for the microscopic scale, we build initial data $\rho_0, u_0$ to run system \eqref{eulereq}.
Denoting $X^0\in \mathbb{R}^N$, $V^0\in \mathbb{R}^N$ initial positions and velocities of the agents, we define an initial distribution probability considering the sum of $N$ Gaussian functions
\begin{equation}
\label{eq:f0_per_macro}
f^0(x,v)= \sum_{i=1}^{N}\frac{1}{2 \pi \sigma_x  \sigma_v} e^{-\frac{(x-X^0_i)^2}{2\sigma_x^2}-\frac{(v-V^0_i)^2}{2\sigma_v^2}} 
\end{equation}
with $\sigma_x=0.1$, $\sigma_v=0.1$. Here $X^0_i , V^0_i \in \mathbb{R}$ denotes the position and the velocity of each agent $i$, respectively.
Recalling the definition of $\nu_t^0$ ad $\nu_t^1$ in \eqref{nu0}-\eqref{nu1}, initial data at the macroscopic scale are given as 
\begin{equation} 
\label{mom_data}
\left\{
\begin{array}{l}
\rho_0=\nu^0_0,\\
Q_{0}= \nu^1_0
\end{array}
\right.
\end{equation}
where $Q_t:=\rho_t u_t$ denotes the momentum, for any $t \ge 0$.

We run \eqref{eulereq} for both choices of the function $K$ in \eqref{function_K1}-\eqref{function_K2}, aiming at comparing the the influence of the performed choice on the convergence time to the consensus state also at the macroscopic scale. We set $T=40$, $\Delta x=0.4$, $\Delta t=0.02$. In particular, we evaluate the macroscopic average velocity $\bar{u}(t)$, computed as the integral of the velocity field over the spatial domain. Figure \ref{Test_K_micromacro}(d) shows that also at this scale the evolution depends on the choice of $K$, with faster convergence to the consensus equilibrium corresponding to the linear one.

\begin{figure}[!h]
\centering
\subfigure{\includegraphics[scale=0.55]{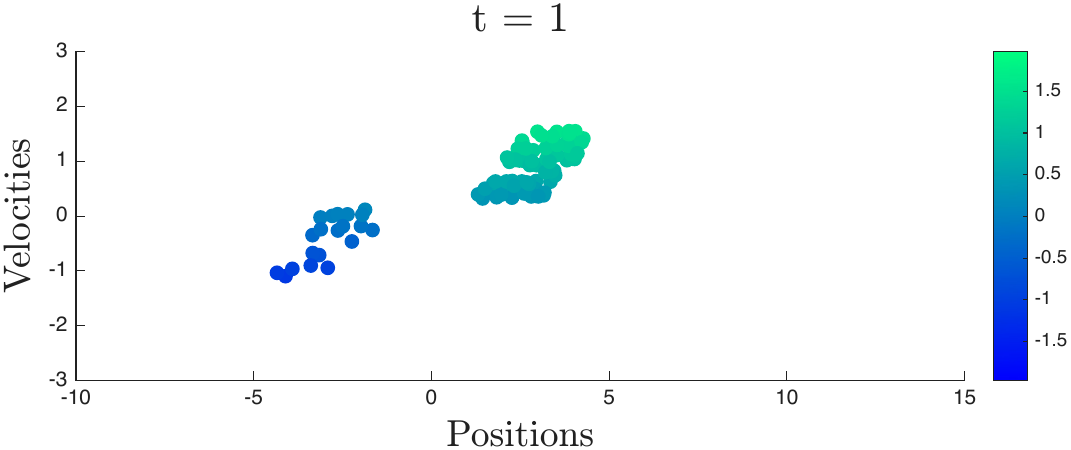}}\\
\subfigure{\includegraphics[scale=0.55]{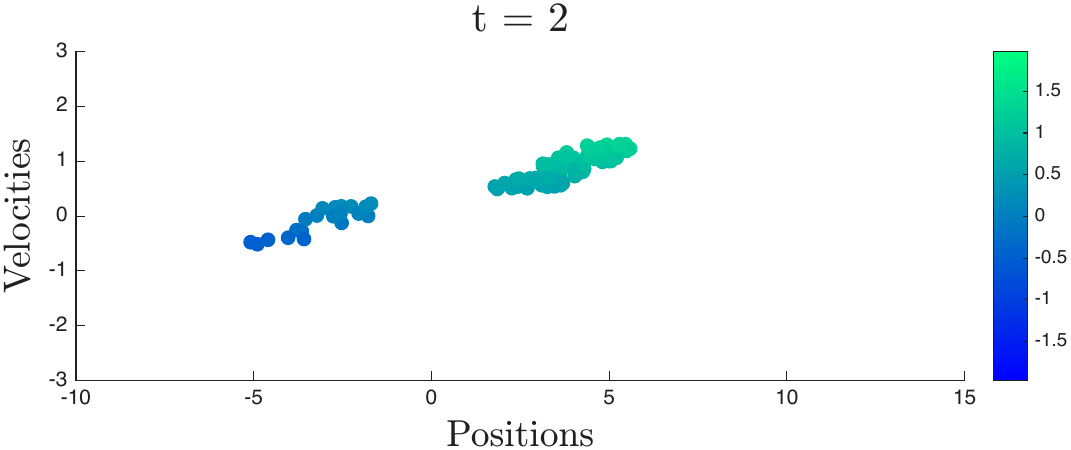}}
\subfigure{\includegraphics[scale=0.55]{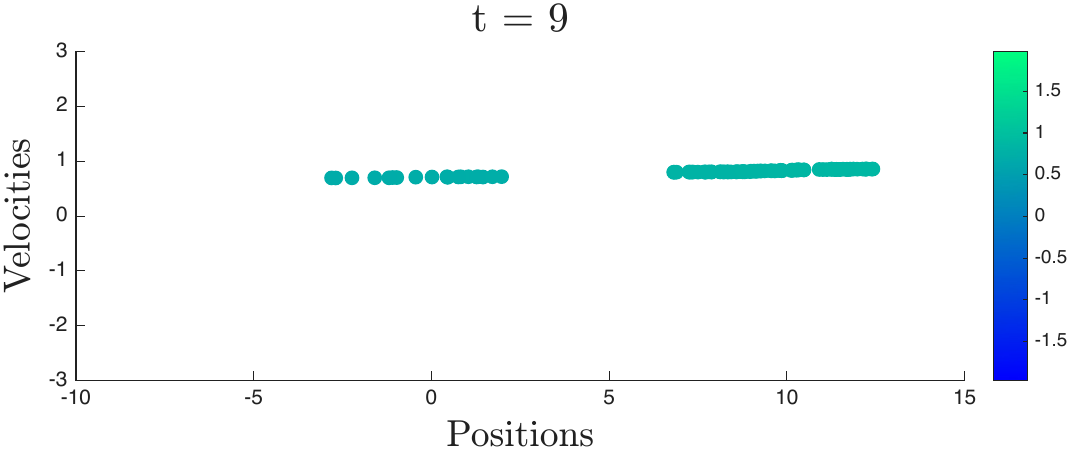}}\\
 \caption{ Test 3: screenshots of numerical simulation of \eqref{sec2:CSmodel_pij} with $K=K_1$ as in \eqref{function_K1}.   }
\label{Test2gruppi_micro_comparison}
\end{figure}

\paragraph{\textit{Test 3: Evolution of the microscopic model for a linear $K$}}
We conclude our analysis on the microscopic scale in 1D showing some screenshots of the dynamics of agents assuming $K=K_1$ in \eqref{function_K1}, see Figure \ref{Test2gruppi_micro_comparison}. 
We remark that initial condition in Test 1 and its evolution in Test 3 represents the reference microscopic dynamics for the comparison with different scales performed in Section \ref{comparison_microkinmacro}.

\newpage
\vspace{10cm}

\subsection{Comparison of the microscopic, kinetic and macroscopic models}
\label{comparison_microkinmacro}

Starting from the results obtained at the microscopic scale, we here focus on the kinetic and macroscopic scales. 
We recall that the derivation of the hydrodynamic limit \eqref{eulereq} is obtained under a monokinetic assumption on the initial data of \eqref{sec2:Vlasov}. 
Indeed, far from the monokinetic case, no rigorous proofs of convergence have been established.

We approximate the solution to \eqref{sec2:Vlasov} in the space-velocity domain $\left[ -15, 20\right] \times \left[ -5, 5\right] $. 
The simulation runs over the time interval $\left[ 0, T\right]=[0, 3]$ with $\Delta t$ chosen according to the CFL condition \eqref{CFL_kinetic}.

\paragraph{\textit{Test 4: Monokinetic initial condition}}
We recall that analytical results obtained in \cite{BPR24} ensure a complete correspondence between Vlasov moments and Euler solutions, in case of monokinetic initial data, namely 
\begin{equation}
    f^{0} = \rho^{0} (x) \delta (v-u^{0}(x)).
\end{equation}
This is a strong \textit{ansatz}, since the solution is completely described by its hydrodynamic fields and no further closure assumptions are required. 
In Test 4 we simulate a monokinetic-like initial data scenario.  We approximate the initial distribution choosing a Gaussian function with a small value of $\sigma_v$, namely:
\begin{equation}
\label{eq:f0_monokin}
f^0(x,v)= \frac{1}{2 \pi \sigma_x  \sigma_v} e^{-\frac{(x-x_0)^2}{2\sigma_x^2}-\frac{(v-v_0)^2}{2\sigma_v^2}} 
\end{equation}
with $x_0=-2$, $v_0=1.5$, $\sigma_x=2$, $\sigma_v=\sqrt{0.001}$. 
The solution is approximated running the scheme detailed in section \ref{sec:numscheme_kinetic} setting $\Delta x = 0.05$ $ \Delta v = 0.005$.
Following the approach in \cite{MNP24}, we compare the moments of the solution to Vlasov with the solution of the Euler system.
We run the macroscopic system \eqref{eulereq} starting with 
\begin{equation} 
\label{mom_data}
\left\{
\begin{array}{l}
\rho_0=\nu^0_0,\\
Q_{0}= \nu^1_0.
\end{array}
\right.
\end{equation}
with $\nu^0_0$, $\nu^1_0$ defined as in \eqref{nu0}-\eqref{nu1}.

\begin{figure}[!h]
\centering
\subfigure{\includegraphics[scale=0.32]{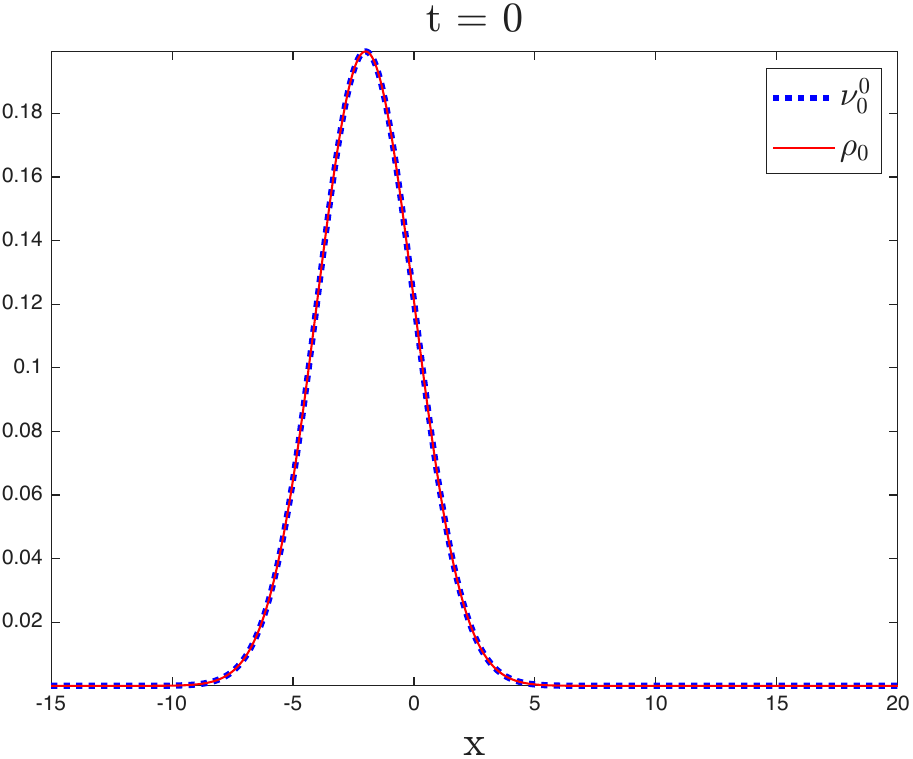}}
\subfigure{\includegraphics[scale=0.32]{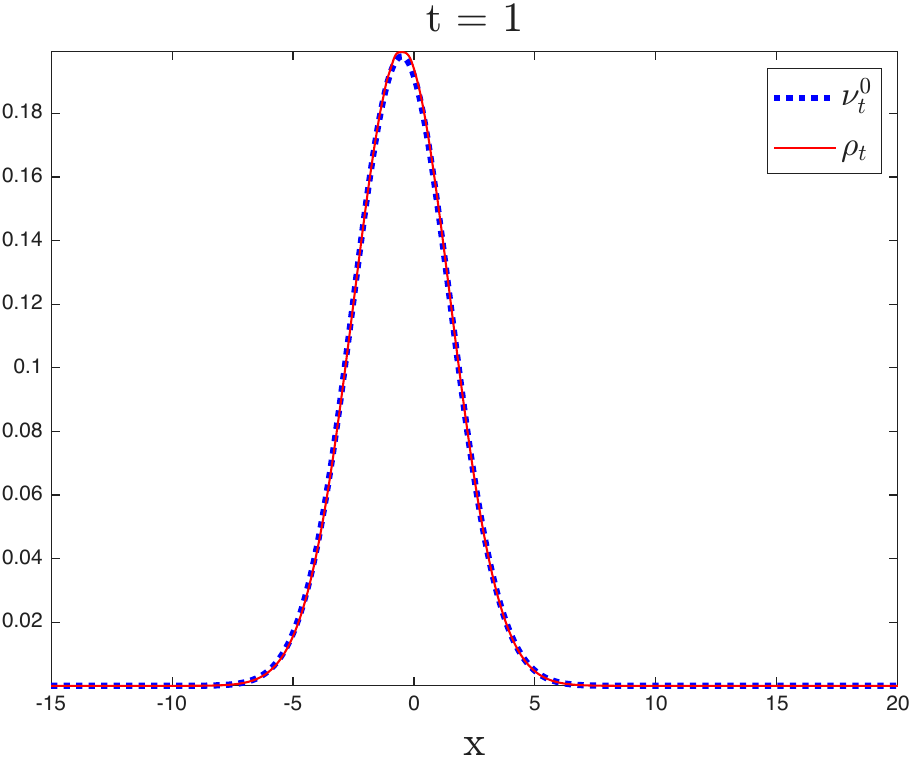}}
\subfigure{\includegraphics[scale=0.32]{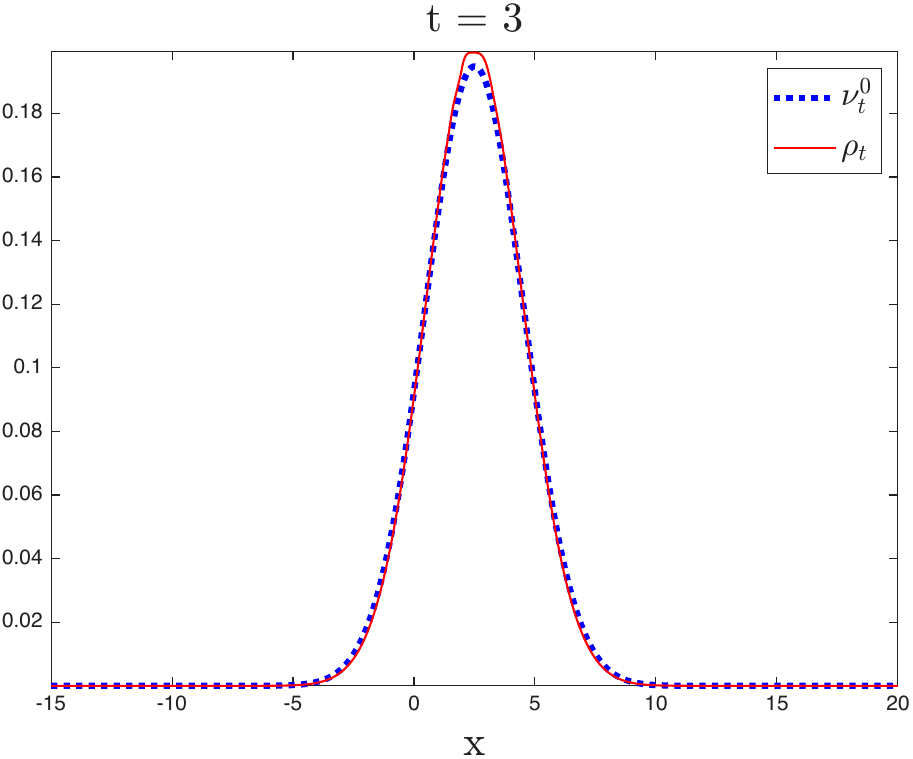}}\\
\subfigure{\includegraphics[scale=0.32]{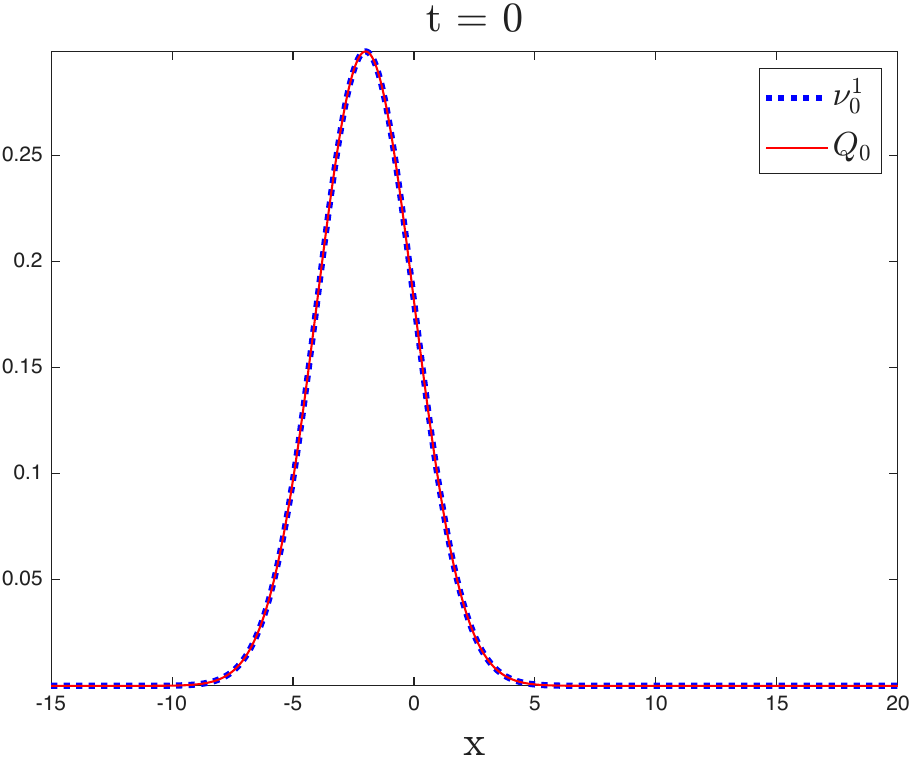}}
\subfigure{\includegraphics[scale=0.32]{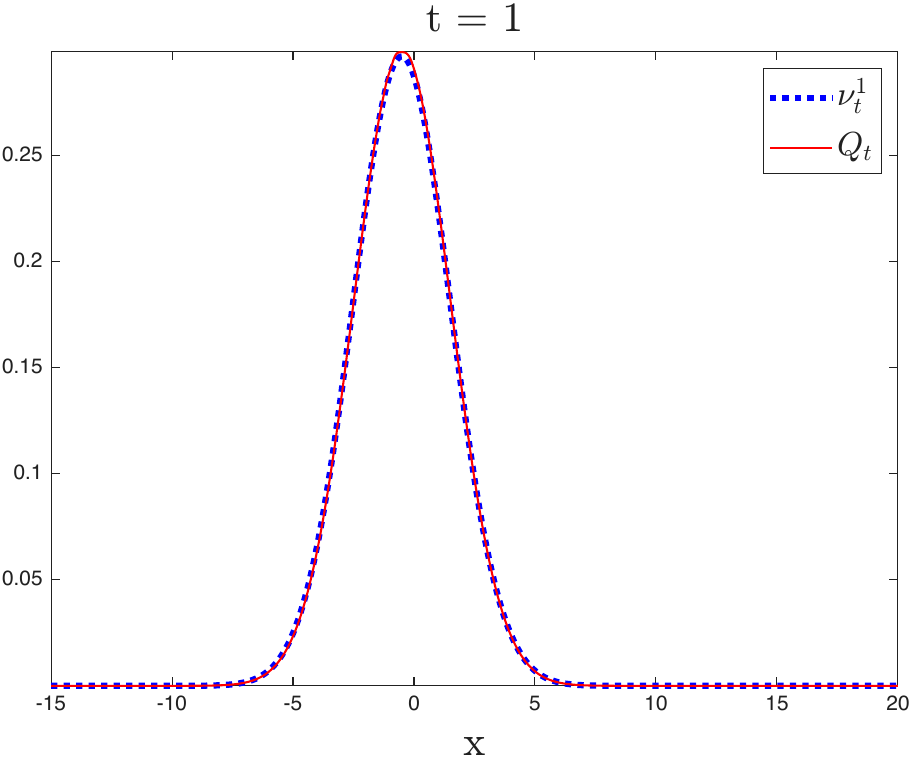}}
\subfigure{\includegraphics[scale=0.32]{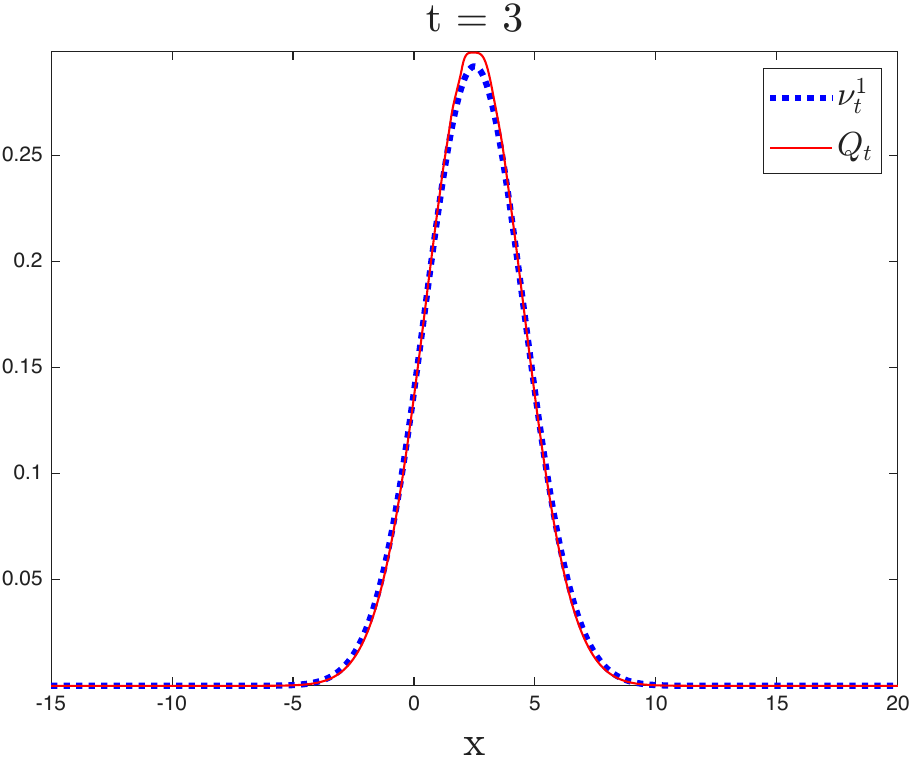}}
 \caption{ Test 4: approximation of monokinetic initial data. Comparison of $\nu_t^0$ and $\rho_t$ (first line) and $\nu_t^1$ and $Q_t$ (second line) at different time steps.}
\label{Test_monokinetic}
\end{figure}

\begin{figure}[!h]
\centering
 \subfigure{\includegraphics[scale=0.4]{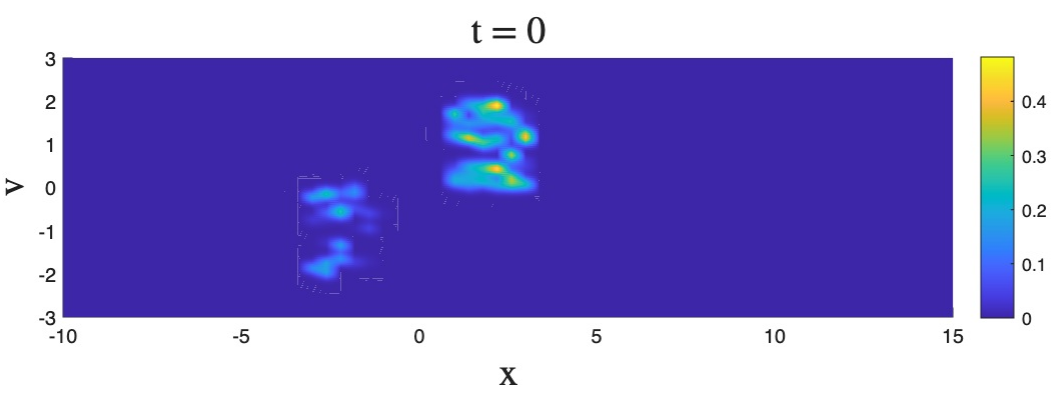}}
\subfigure{\includegraphics[scale=0.4]{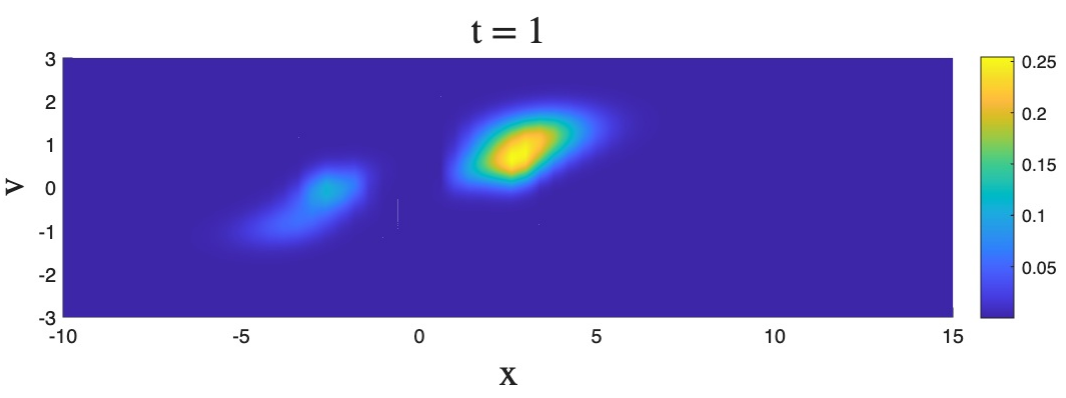}}\\
\subfigure{\includegraphics[scale=0.4]{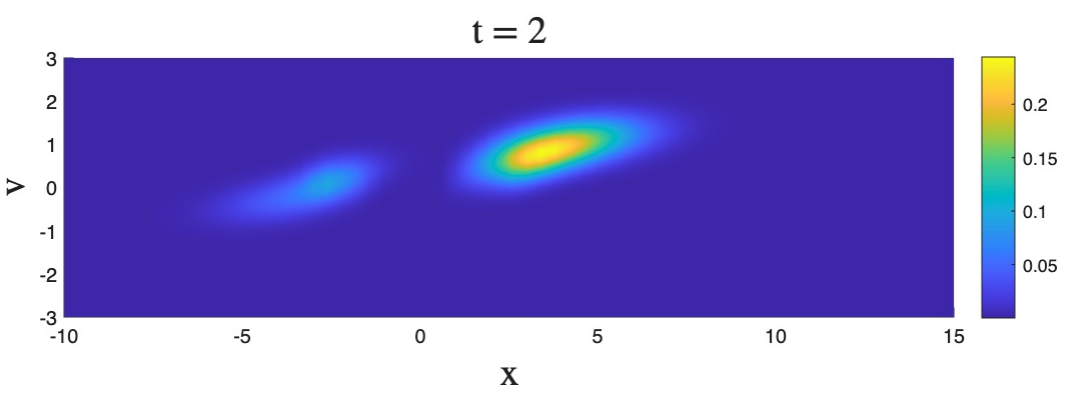}}
\subfigure{\includegraphics[scale=0.4]{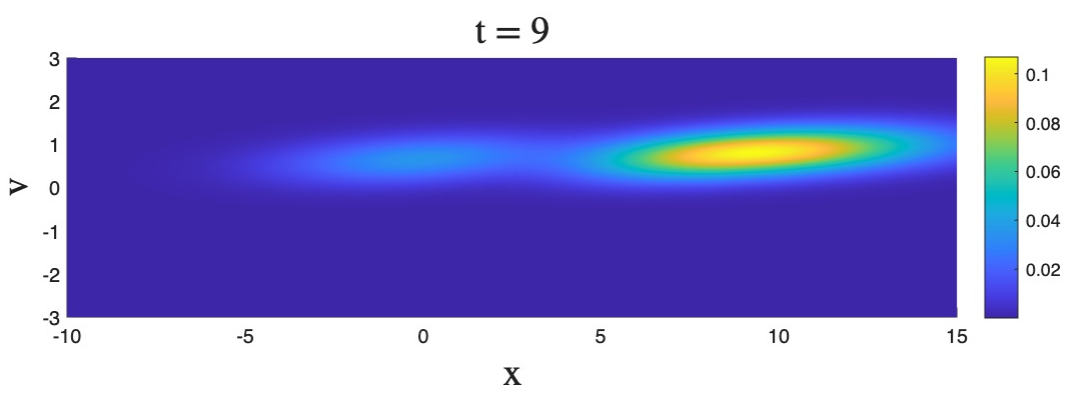}}
 \caption{ Test 5: screenshots of a numerical simulations of \eqref{sec2:Vlasov} in the phase space. }
\label{Test2gruppi_kin_phasespace}
\end{figure}

\begin{figure}[!h]
\centering
\subfigure{\includegraphics[scale=0.32]{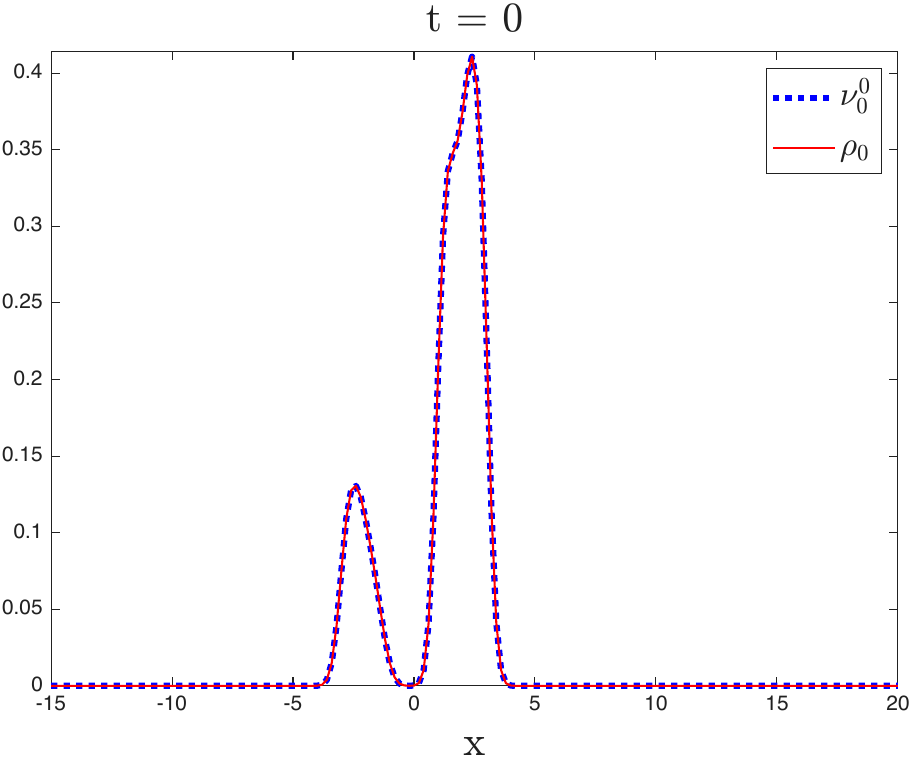}}
\subfigure{\includegraphics[scale=0.32]{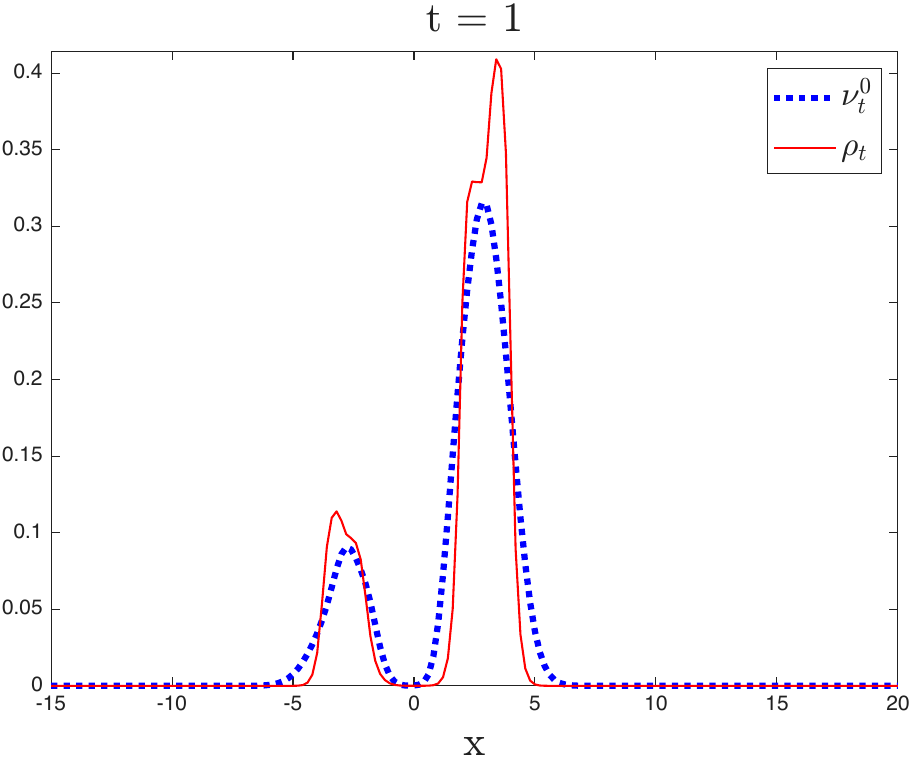}}
\subfigure{\includegraphics[scale=0.32]{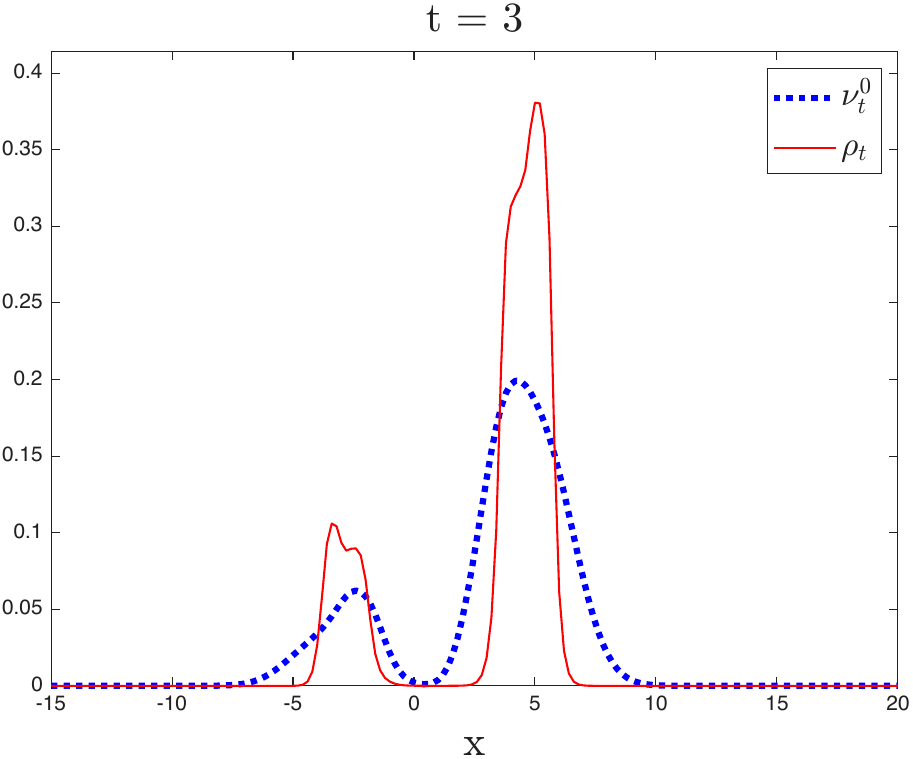}}\\
\subfigure{\includegraphics[scale=0.32]{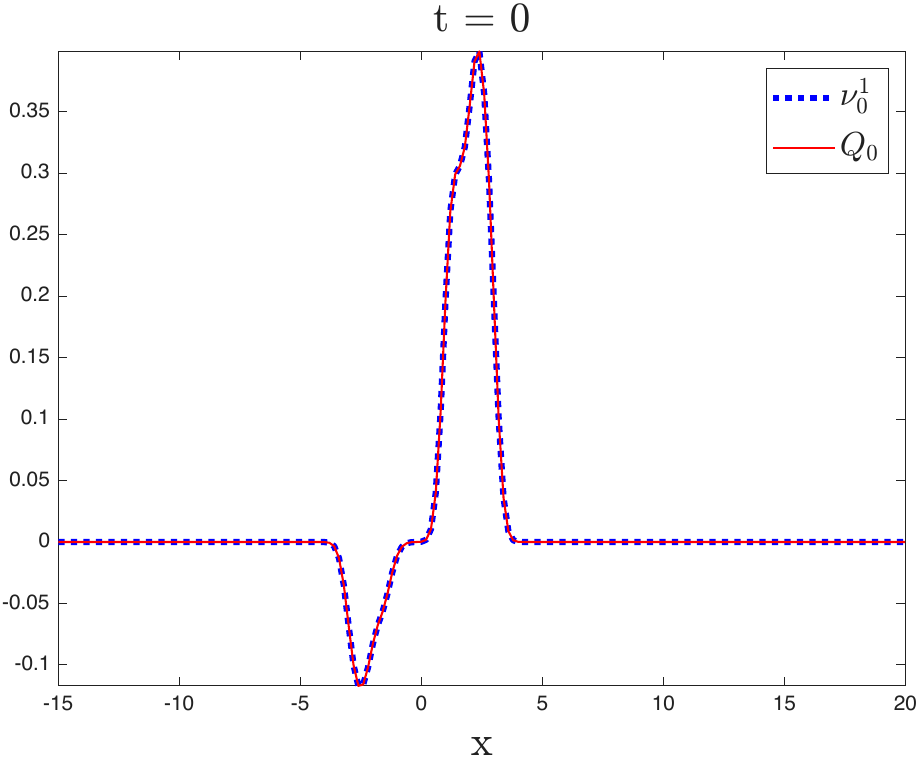}}
\subfigure{\includegraphics[scale=0.32]{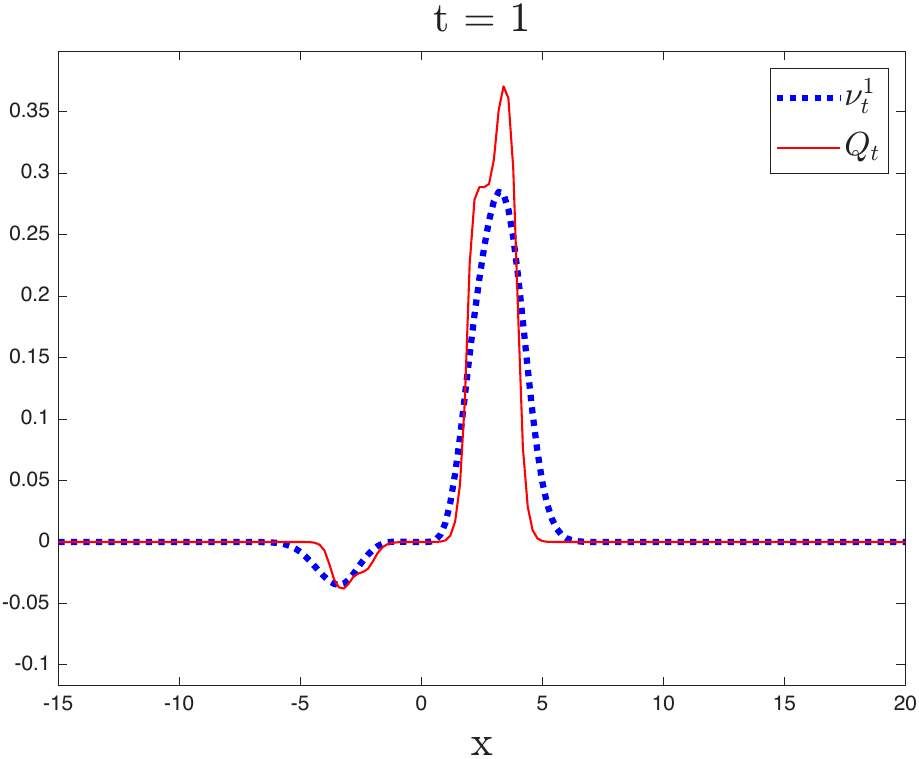}}
\subfigure{\includegraphics[scale=0.32]{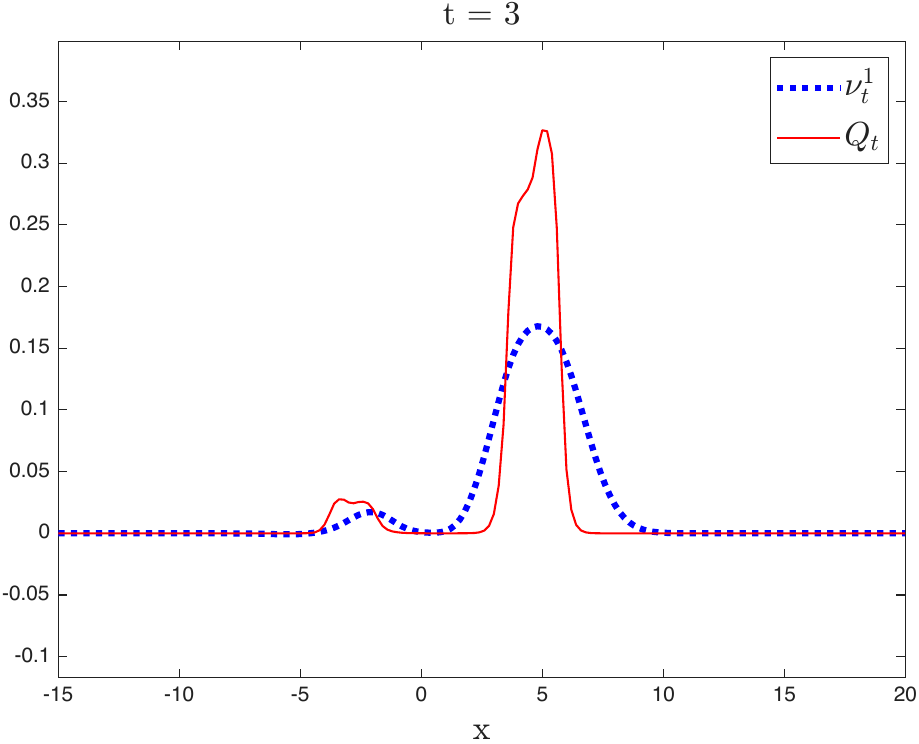}}
 \caption{ Test 5: Comparison of $\nu_t^0$ and $\rho_t$ (first line) and $\nu_t^1$ and $Q_t$ (second line) at different time steps.}
\label{Test_NONmonokin_kinmacro}
\end{figure}

Figure \ref{Test_monokinetic} shows the comparison between the scales at two different time instants. As expected, we observe a good agreement between zeroth and first order moments of the solution to \eqref{sec2:Vlasov}, and density and momentum solutions to \eqref{eulereq}. The small discrepancy between the profiles is due to the approximation of monokineticity assumption.
\paragraph{\textit{Test 5: Non-monokinetic initial condition}}
In Test 5 we consider a general initial data, to perform a comparison far from the monokinetic assumption. We run \eqref{sec2:Vlasov} setting $\Delta x = 0.2$ $ \Delta v = 0.025$, starting with the initial distribution defined in \eqref{eq:f0_per_macro}.
Screenshots of the distribution at different times are reported in the phase space, see Figure \ref{Test2gruppi_kin_phasespace}.
Hence we run the macroscopic model \eqref{eulereq} with initial condition given by \eqref{mom_data}. 
Figure \ref{Test_NONmonokin_kinmacro} shows the comparison between kinetic and macroscopic scale, at different time instant. The non-monokinetic structure of the initial data seems to play a crucial role. Indeed the agreement between kinetic moments and macroscopic solution is no longer observed. Numerical results seem to suggest that the comparison is still valid only for short time, getting worse for longer time.

\newpage

\subsection{Sensitive dependence on initial data: numerical results in 1D and 2D}

\begin{figure}[!h]
\centering
\subfigure[$\varepsilon=-0.04$]{\includegraphics[scale=0.35]{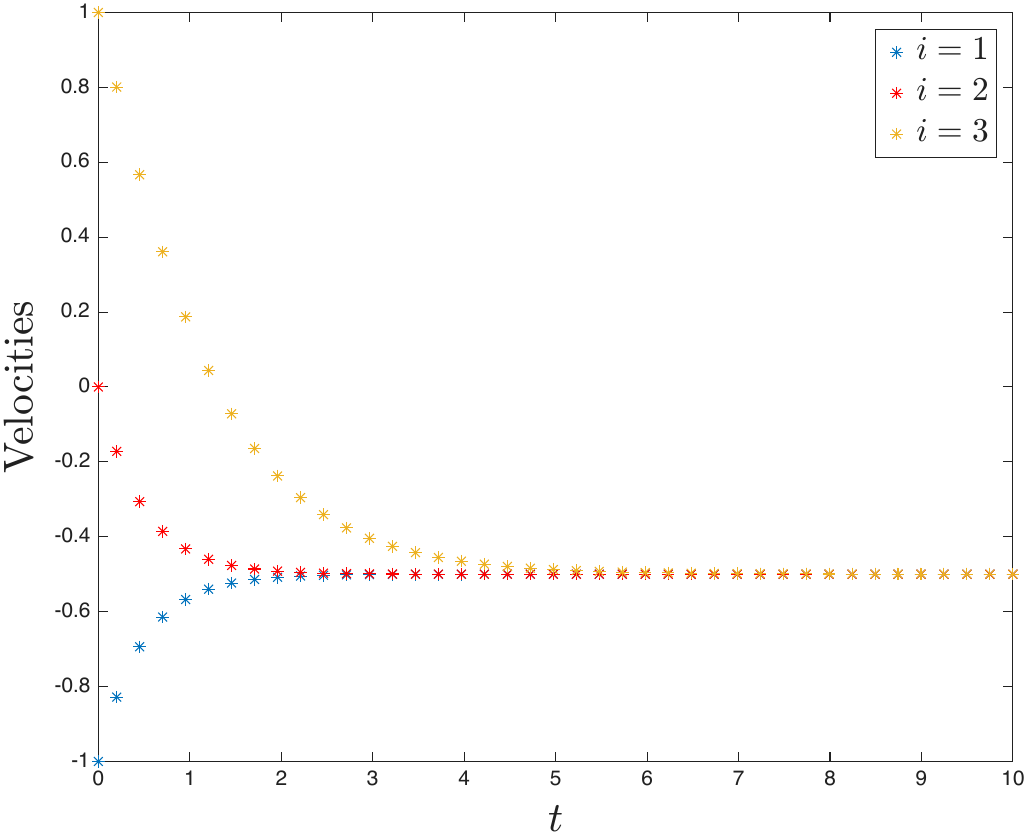}}
\subfigure[$\varepsilon=0.04$]{\includegraphics[scale=0.35]{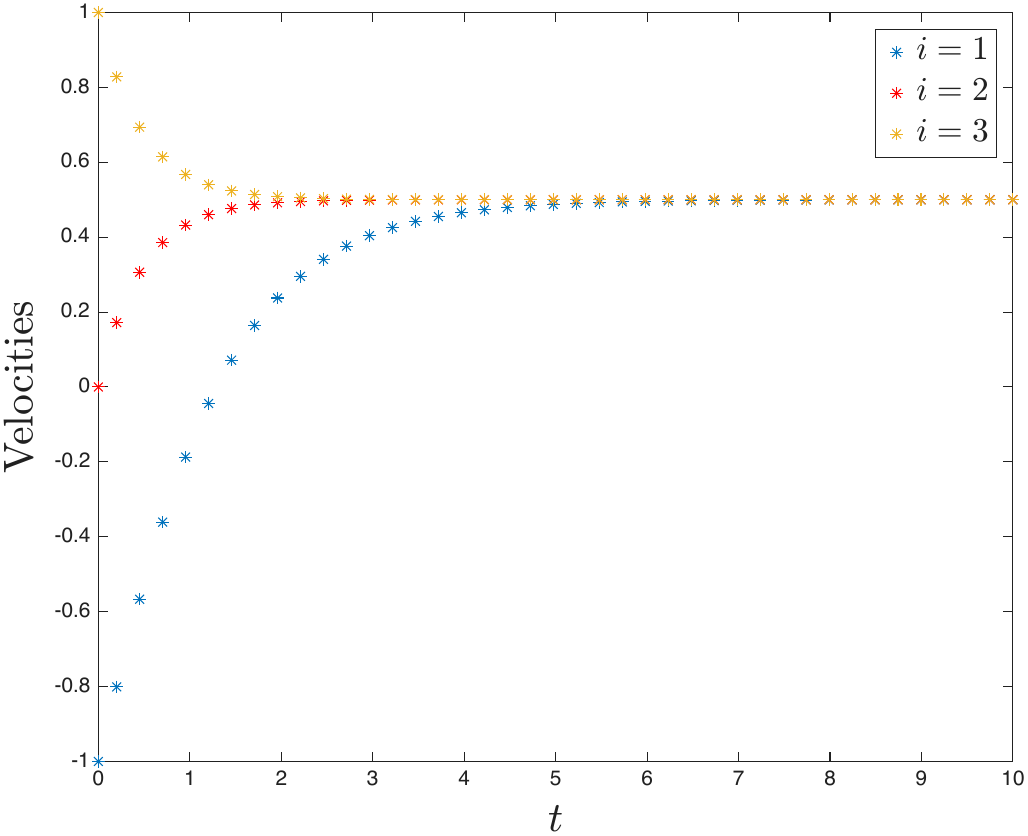}}
 \caption{Test 6: numerical simulation in 1D of \eqref{sec2:CSmodel_pij}. \REV{Here N=3,} $X^0=(-1, \varepsilon, 1)$, $V^0=(-1, 0, 1)$.  }
\label{micro_Rossi}
\end{figure}

\begin{figure}[!h]
\centering
\includegraphics[scale=0.35]{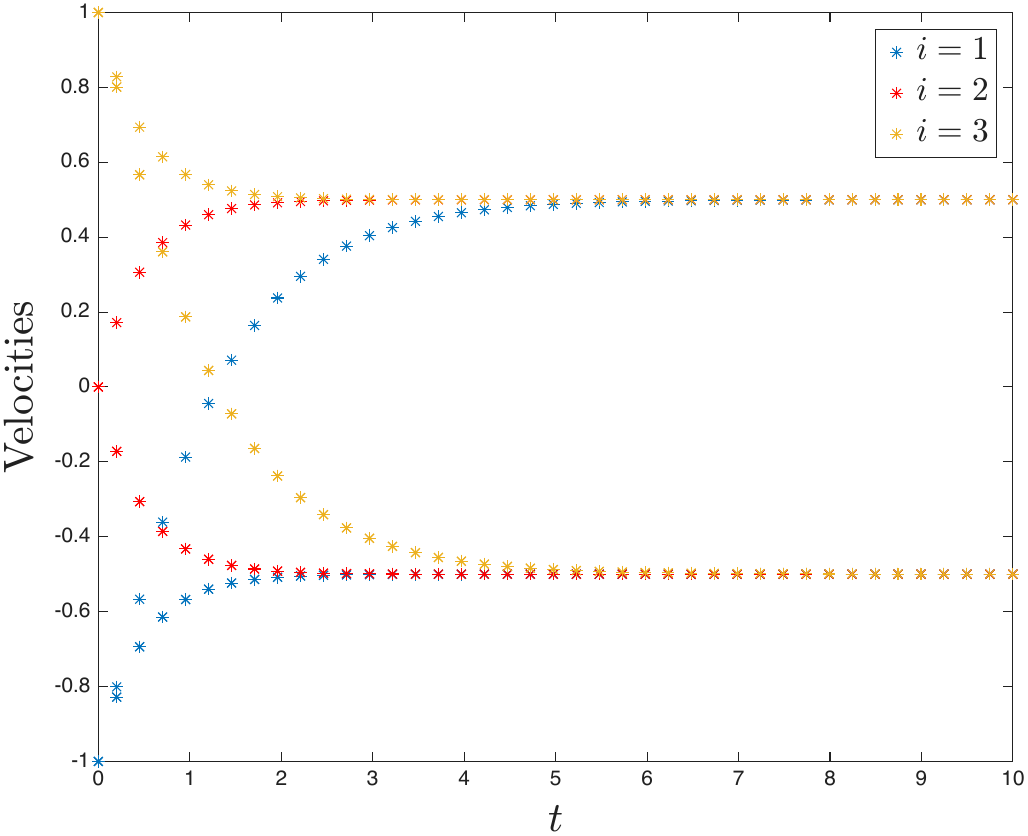}
 \caption{Test 6: numerical simulation in 1D of \eqref{sec2:CSmodel_pij}. Here N=3, $X^0=(-1, 0, 1)$, $V^0=(-1, 0, 1)$, depending on the two \textit{by hand} selections of the first neighbor.}
\label{fig:superposition}
\end{figure}

We conclude our multiscale analysis focusing on the sensitive dependence of the model on initial conditions. 
Indeed, as shown in \cite{BCR21}, there exist pathological initial configurations for which the dynamics is not well-defined. Specifically, this occurs when two agents are equidistant from a third agent, and this equidistance is preserved in the past or future.
\paragraph{\textit{Test 6: Sensitive dependence on initial data in $1$D}}
In \cite{BCR21}, authors prove well-posedeness of \eqref{sec2:CSmodel_pij} for almost all initial data, giving an example of critical initial datum for the system in 1D. 
Initial positions and velocities with $N=3$ agents are given as $X^0=(-1, \varepsilon, 1)$, $V^0=(-1, 0, 1)$, with $\varepsilon \in (-1,1), \varepsilon \neq 0 $. The model coefficients $p_{ij}$ take values $\{ 1/3, 2/3, 1 \}$. Assuming $K(2/3)=1$ and $K(1)=0$, the dynamics varies depending on $\varepsilon$, leading to a discontinuity in $\varepsilon=0$. 
Figure \ref{micro_Rossi} shows the evolution of $V_i$, for $i=1,2,3$, for two different values of the initial position of agent labeled with $i=2$.
We observe that a small variation on the value of $\varepsilon$ reflects in a fairly different dynamics.
In particular, for $\varepsilon = 0$, the solution to system \eqref{sec2:CSmodel_pij} is not unique, since there exists at least two solutions, each being approached by the solution obtained for $\varepsilon \to 0^+$ and $\varepsilon \to 0^-$. Figure \ref{fig:superposition} shows the two solutions obtained running the system assuming $\varepsilon=0$. Since in this case the first nearest neighbor is not unique, we obtain two different solutions, forcing the algorithm to select the particle located at $X^0_1$ or the one at $X^0_3$, leading to two different consensus velocities.
This is the same phenomena observed in standard one particle dynamics with low regularity (less than Cauchy-Lipschitz), see e.g. \cite[Remark IV.3]{LionsPaul}.
It is important to note that — in the one dimensional example considered — the ambiguity in selecting the nearest agent persists over time. Indeed, in this case, two agents are equidistant from a third one, and the equidistance is preserved in the past or future. This situation corresponds to what are defined as irregular points of the \textit{iso-rank manifold} in \cite{BCR21}. Therefore, any selection criterion, such as choosing the agent with the lowest (or highest) index, can only be arbitrary.
It is worth noting that in \cite{BCR21} and \cite{BPR24}, the derivation of the kinetic equation is performed outside this set of initial conditions which
are pathological for the well-posedness of the dynamics and have Lebesgue measure zero.

\begin{figure}[!h]
\centering
\subfigure{\includegraphics[scale=0.3]{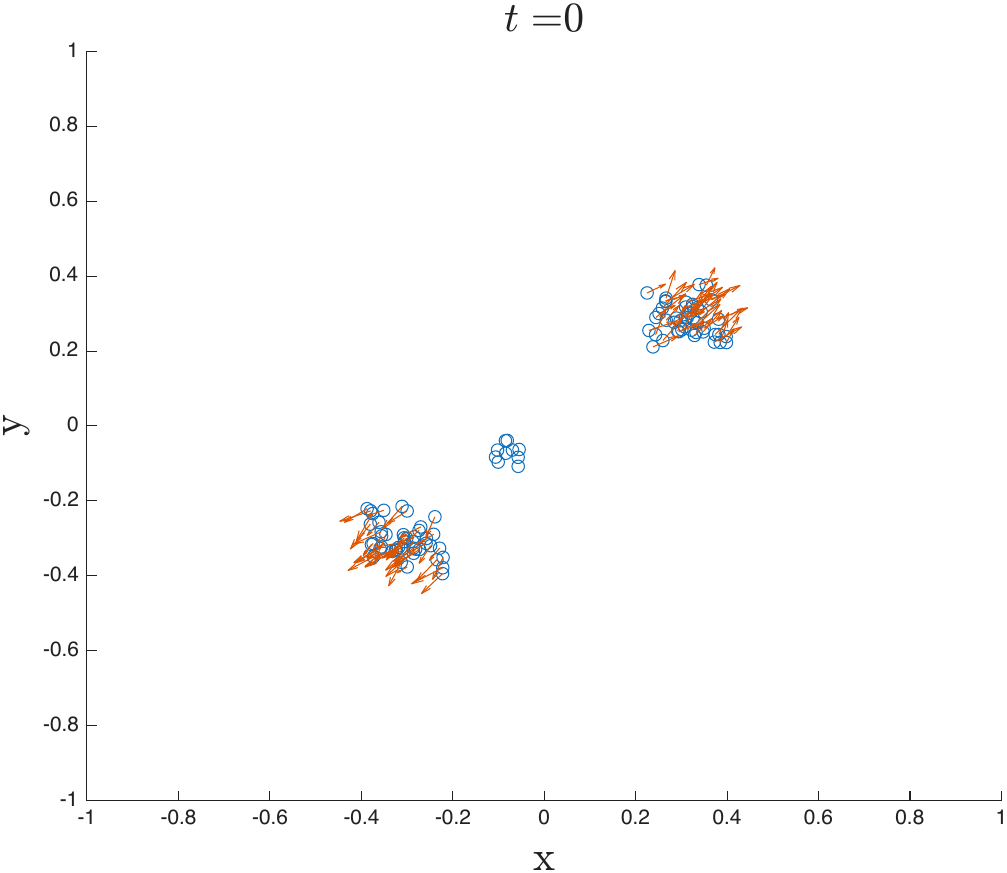}}
\subfigure{\includegraphics[scale=0.3]{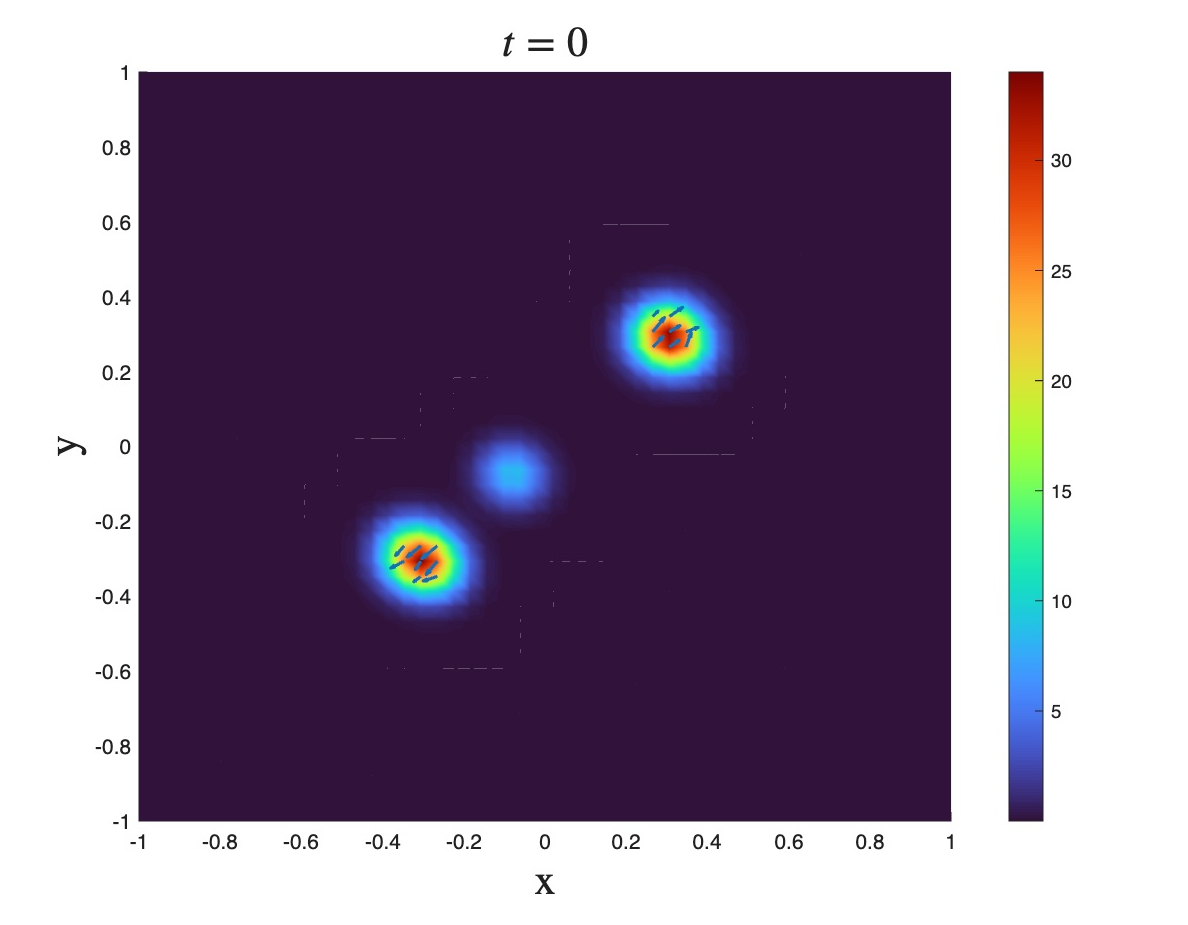}}
 \caption{Test 7: initial conditions for  the microscopic scale and  the macroscopic one. The small group is located in a circle centered in $(-0.08,-0.08)$. }
\label{CI_2D}
\end{figure}

\begin{figure}[!h]
\hspace{-0.8cm}
\subfigure{\includegraphics[scale=0.26]{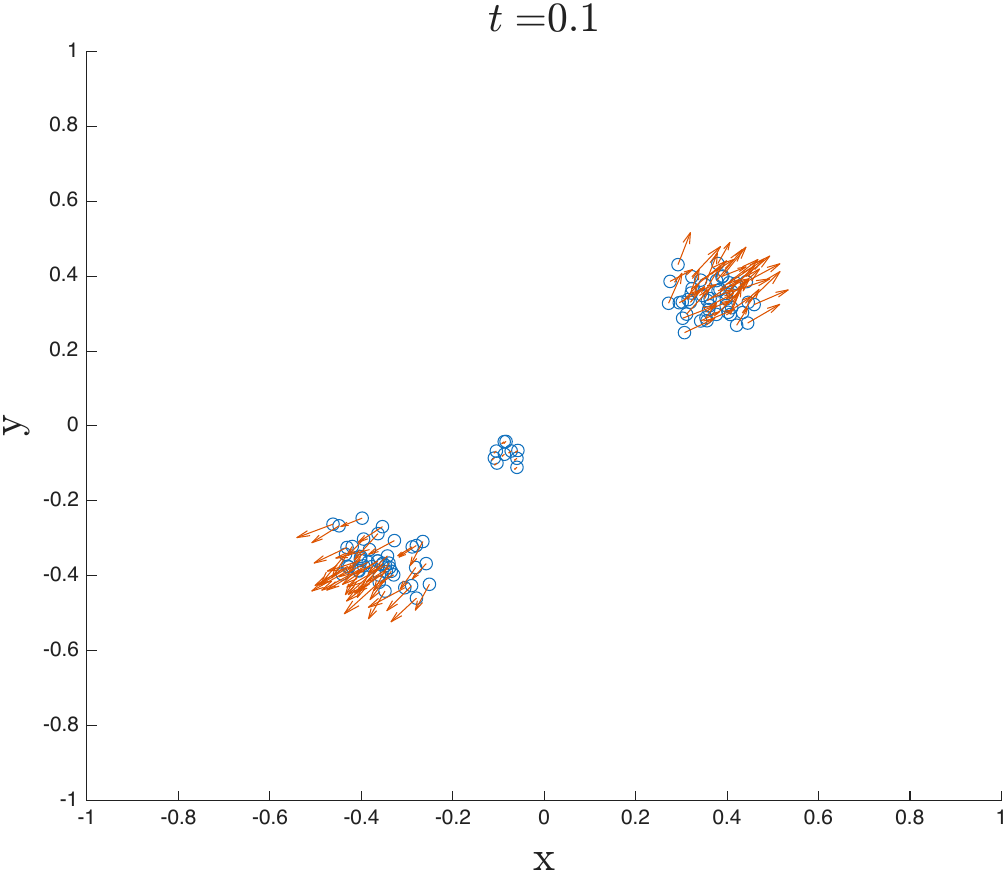}}
\hspace{0.7cm}
\subfigure{\includegraphics[scale=0.26]{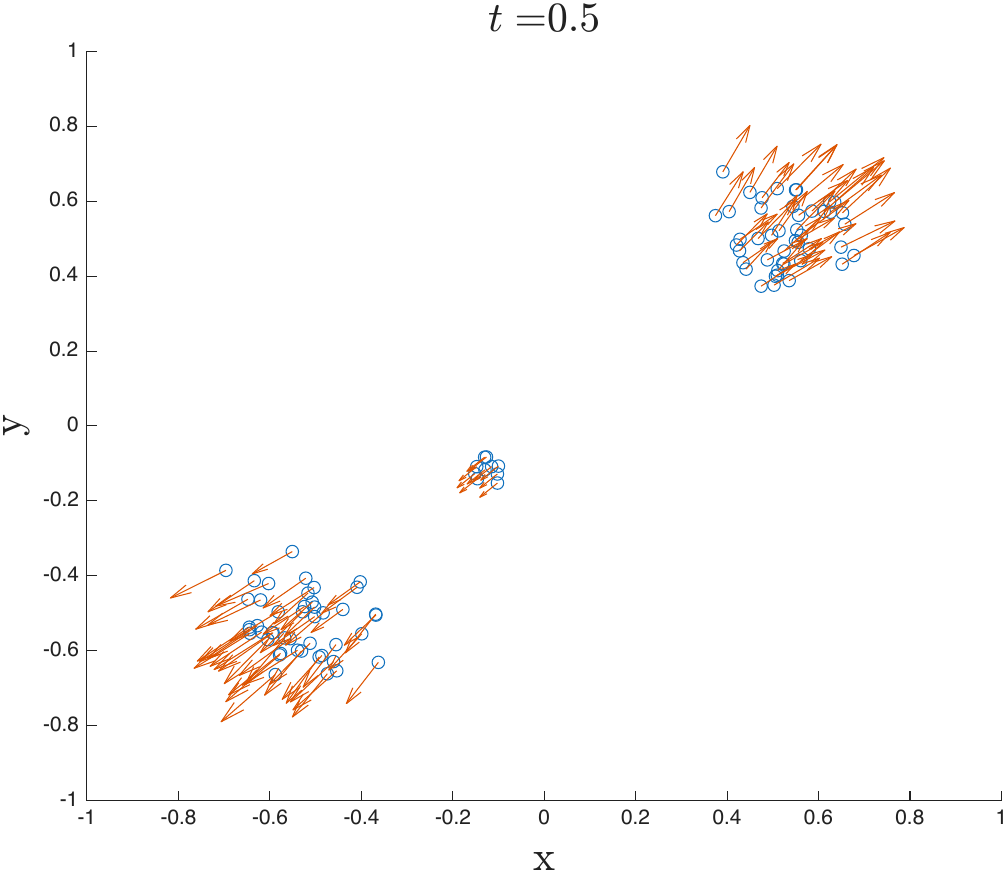}}
\hspace{0.7cm}
\subfigure{\includegraphics[scale=0.26]{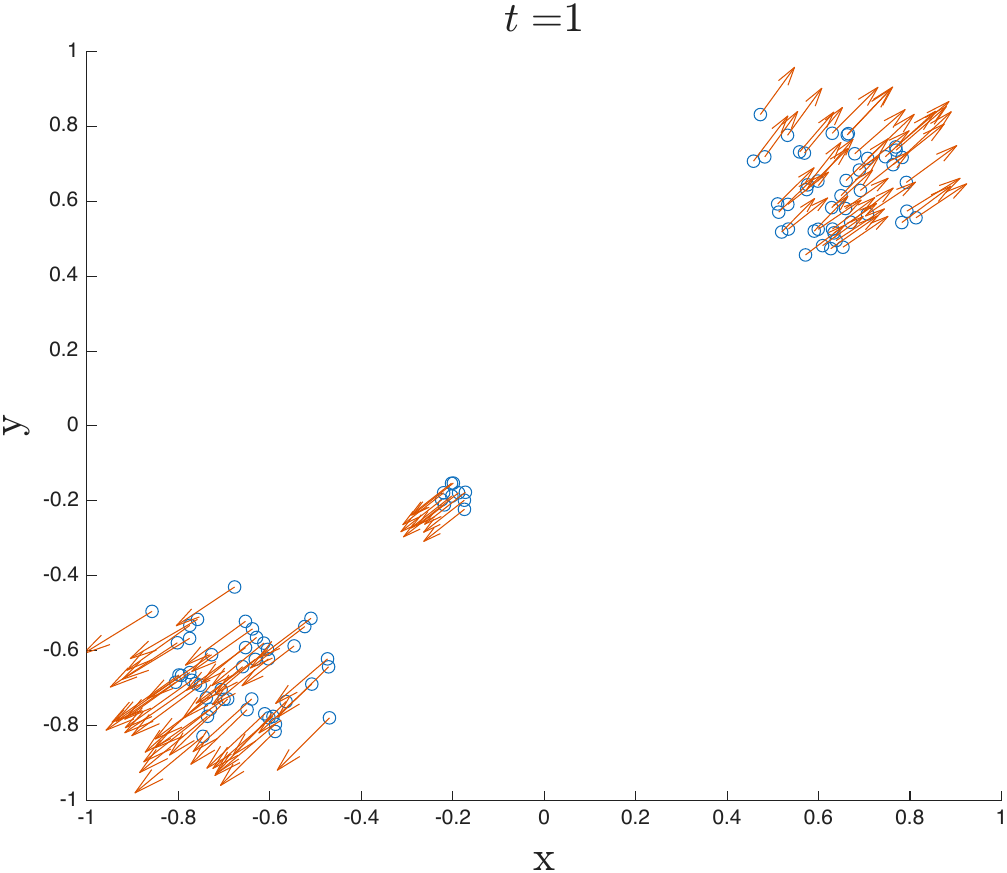}}\\
\subfigure{\includegraphics[scale=0.26]{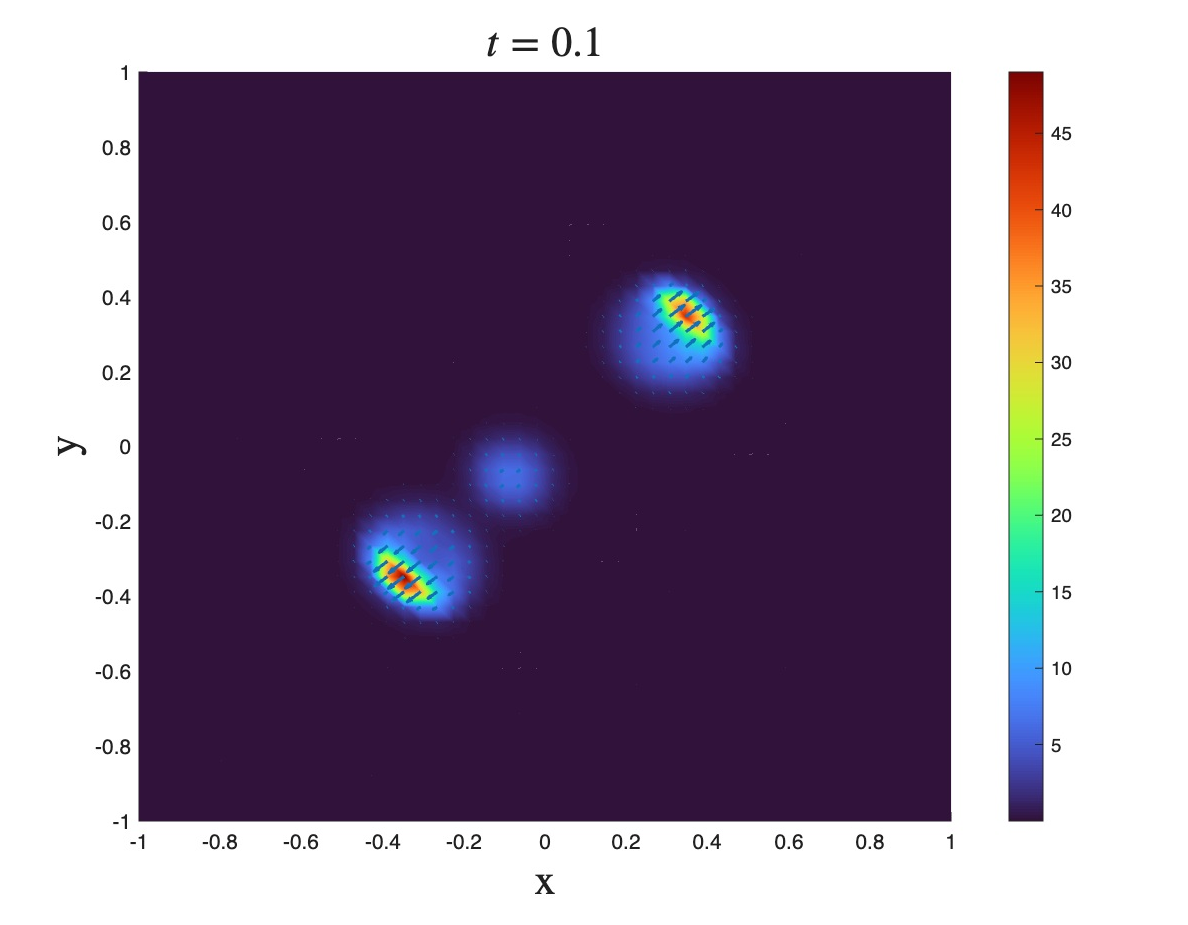}}
\subfigure{\includegraphics[scale=0.26]{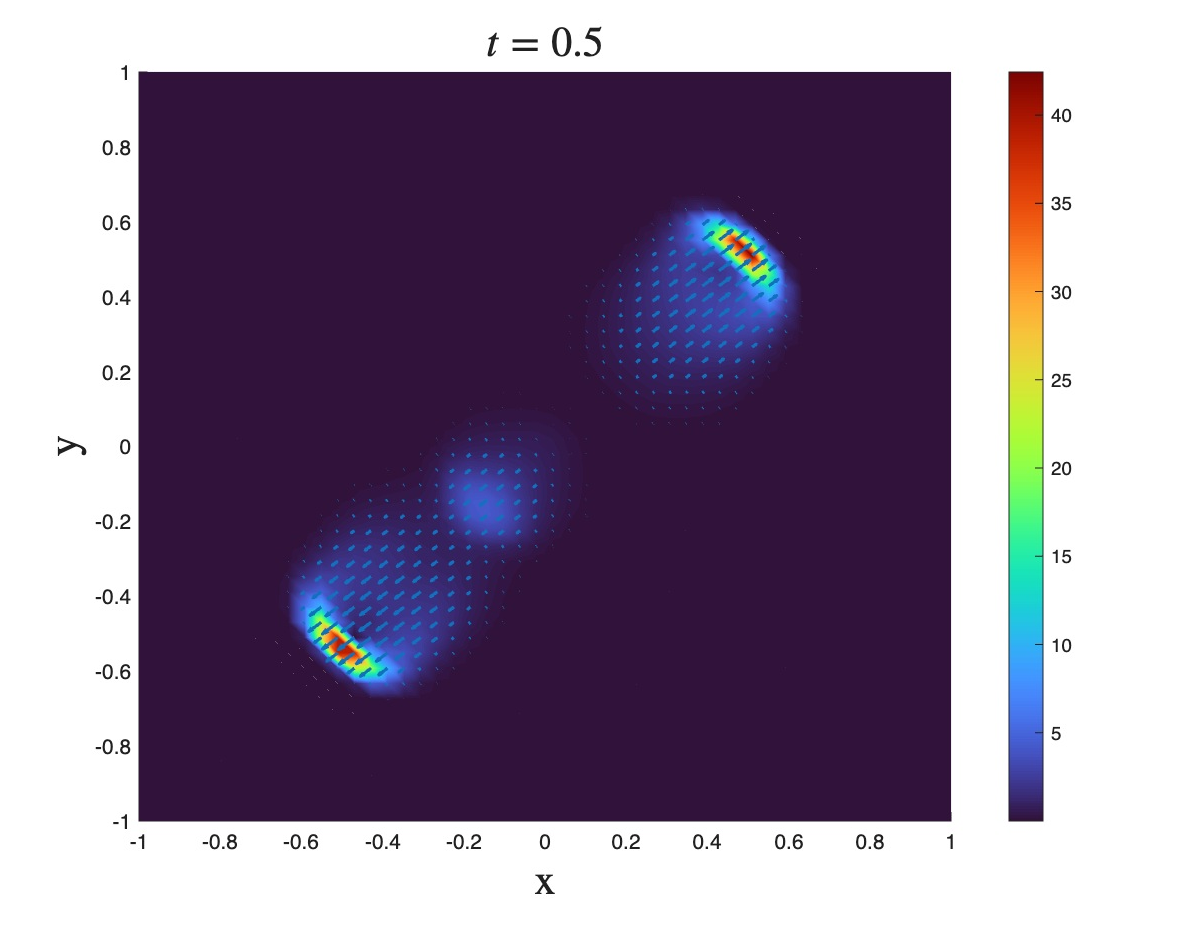}}
\subfigure{\includegraphics[scale=0.26]{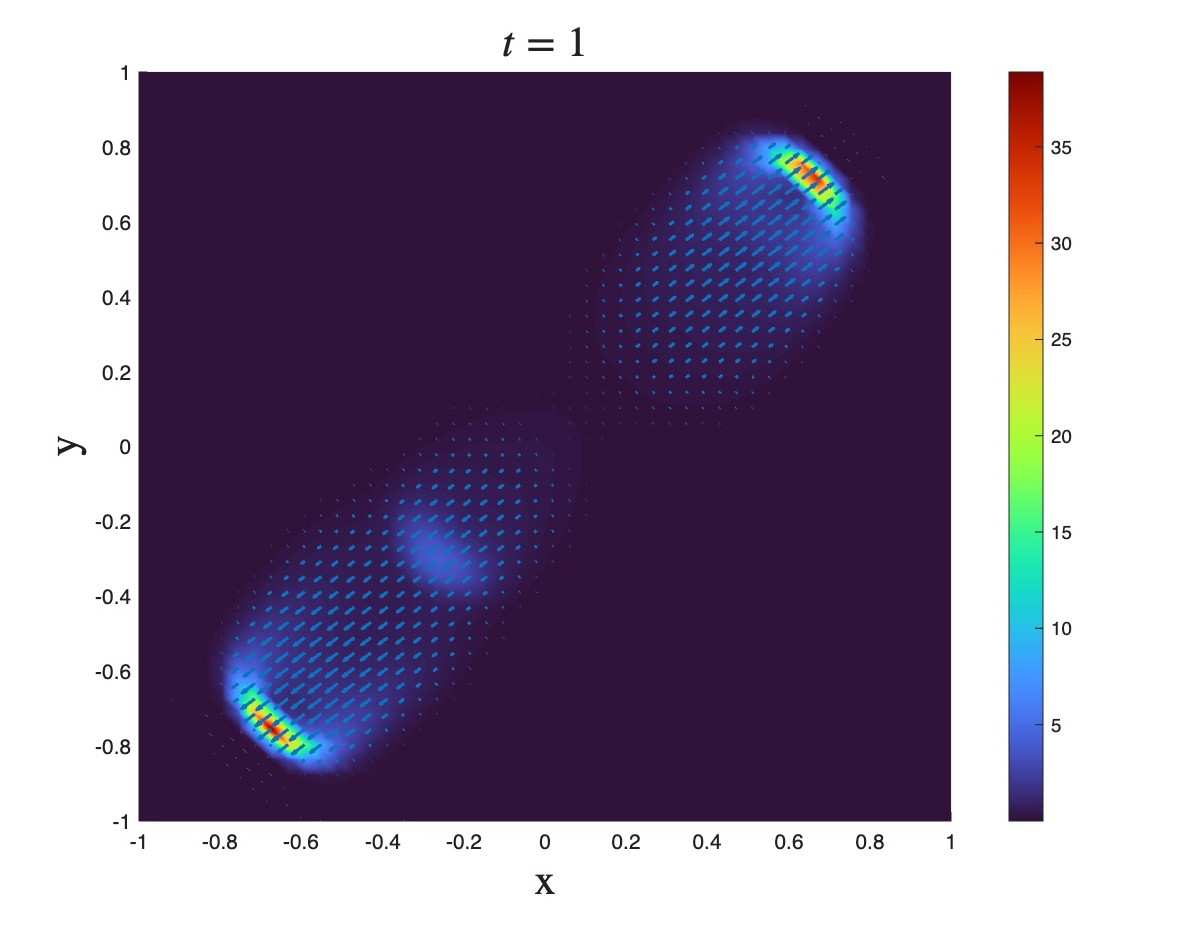}}
 \caption{ Test 7: screenshots of the numerical simulation at (first line) the microscopic scale and (second line) the macroscopic scale.  }
\label{Test_2D}
\end{figure}

\begin{figure}[!h]
\centering
\subfigure{\includegraphics[scale=0.3]{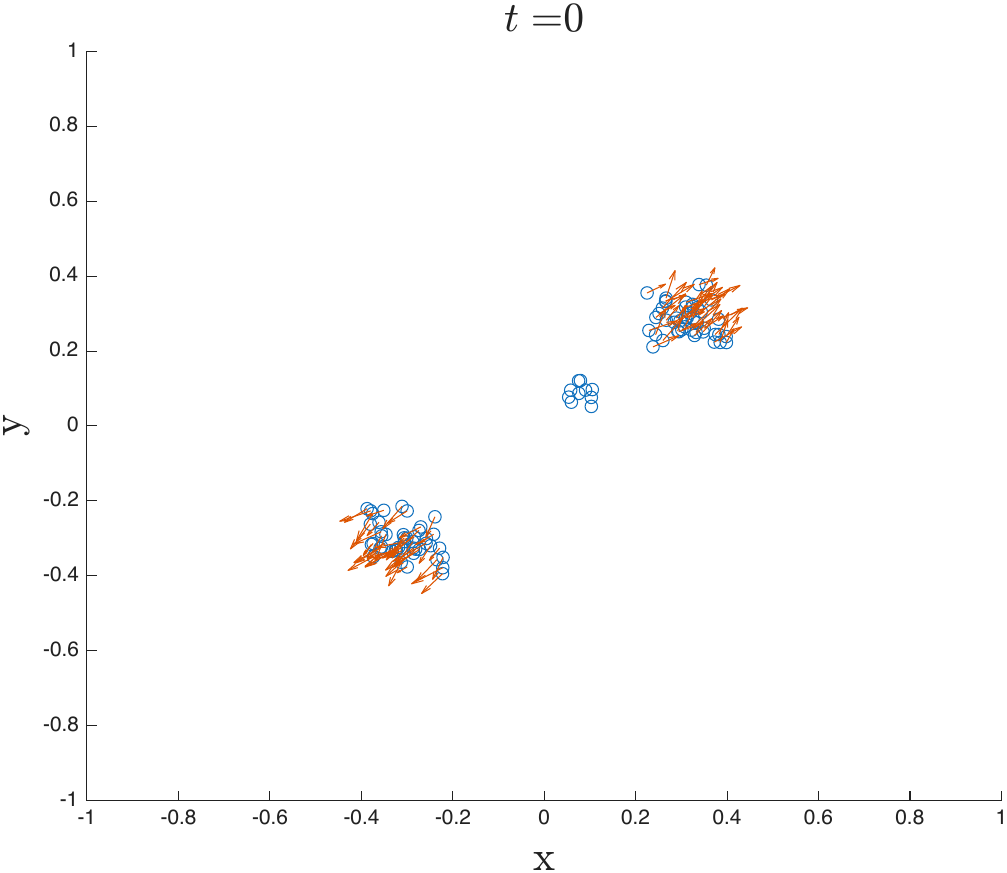}}
\subfigure{\includegraphics[scale=0.3]{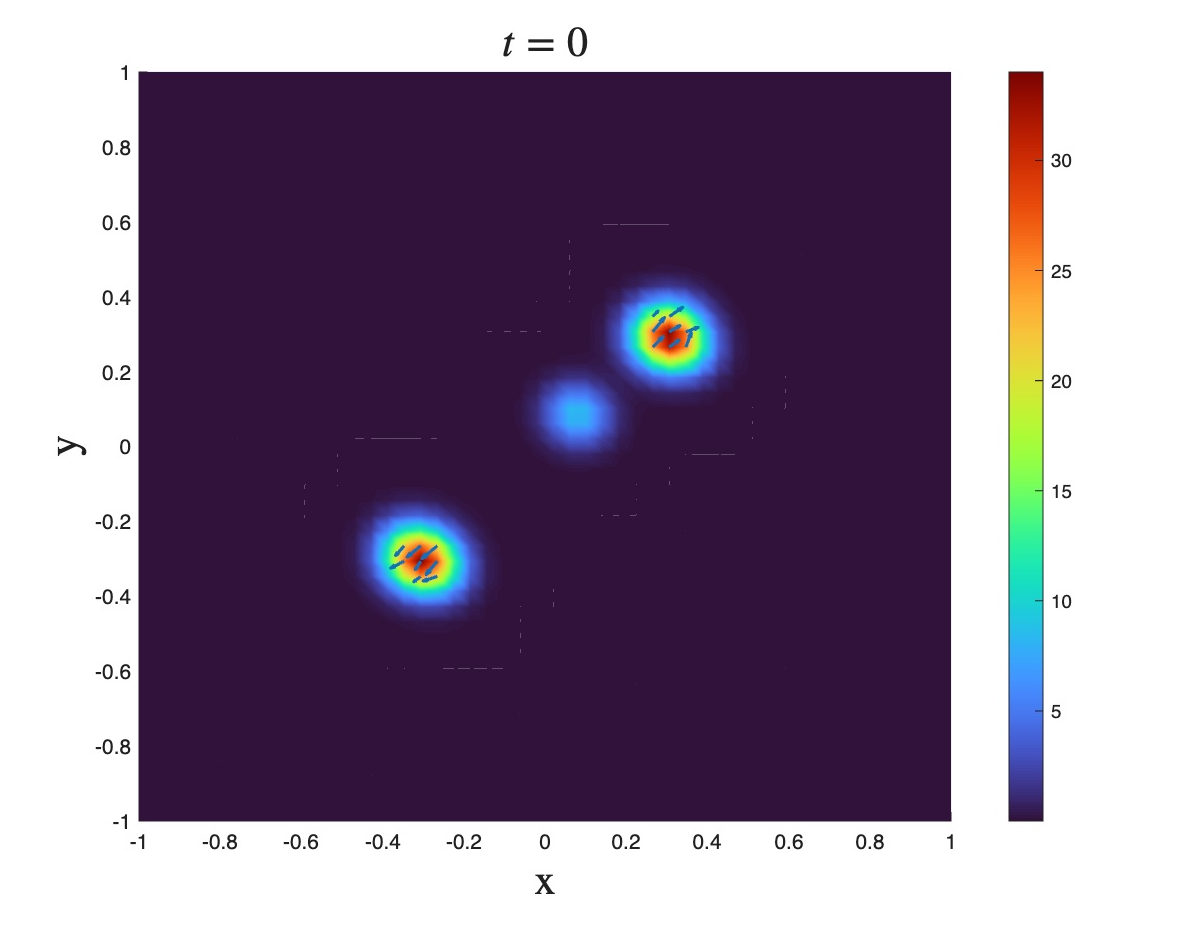}}
 \caption{ Test 8: initial conditions for the microscopic scale and the macroscopic one. The small group is located in a circle centered in $(0.08,0.08)$.   }
\label{CI_2D_positive}
\end{figure}

\begin{figure}[!h]
\hspace{-0.8cm}
\subfigure{\includegraphics[scale=0.26]{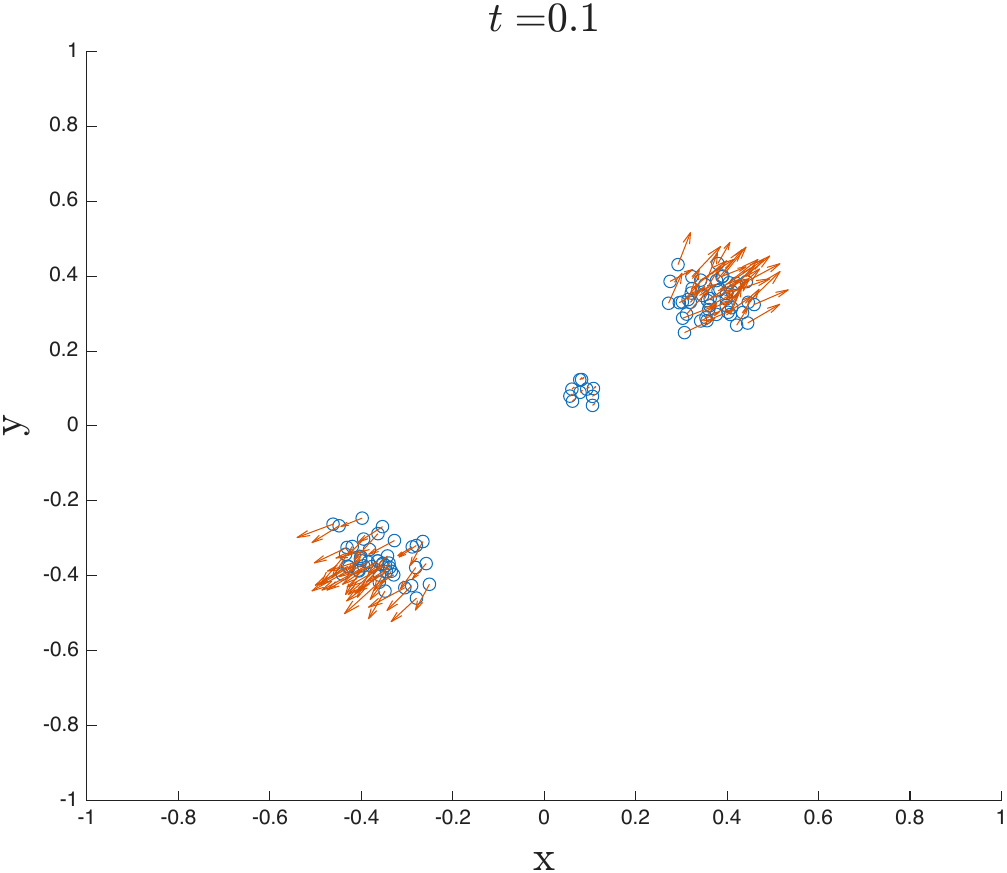}}
\hspace{0.7cm}
\subfigure{\includegraphics[scale=0.26]{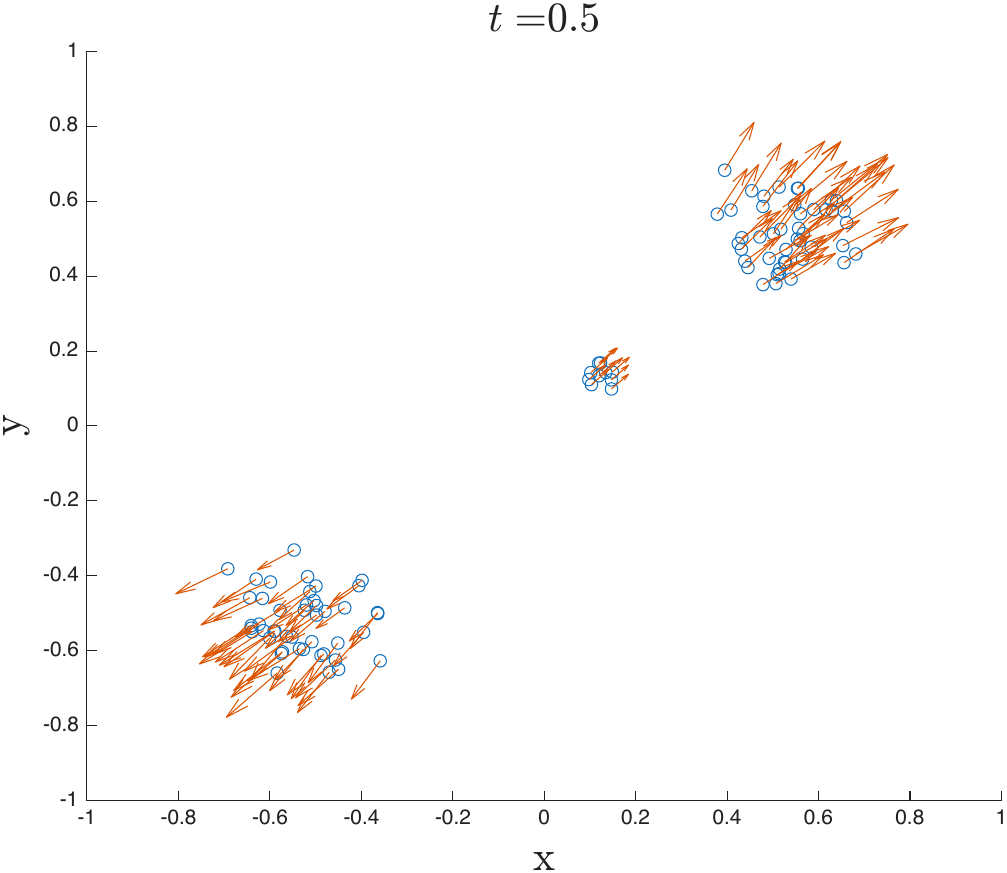}}
\hspace{0.7cm}
\subfigure{\includegraphics[scale=0.26]{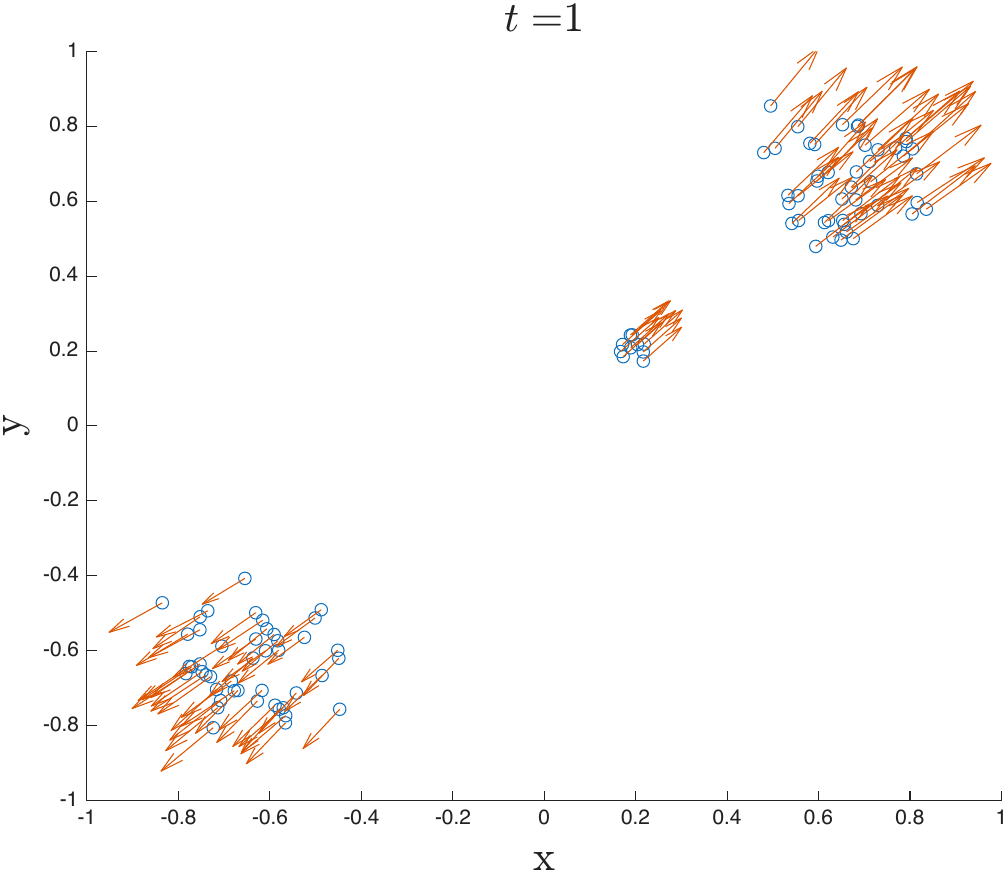}}\\
\subfigure{\includegraphics[scale=0.26]{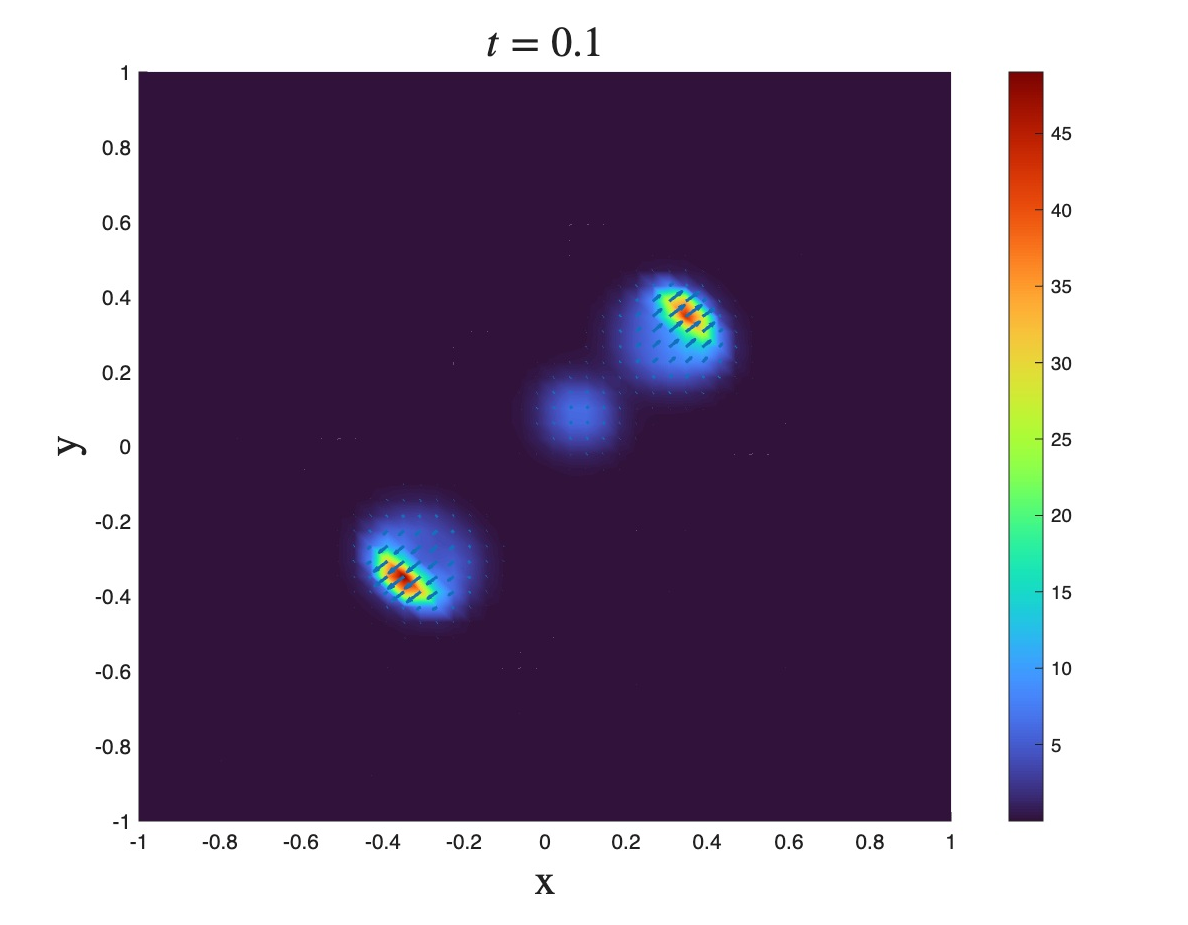}}
\subfigure{\includegraphics[scale=0.26]{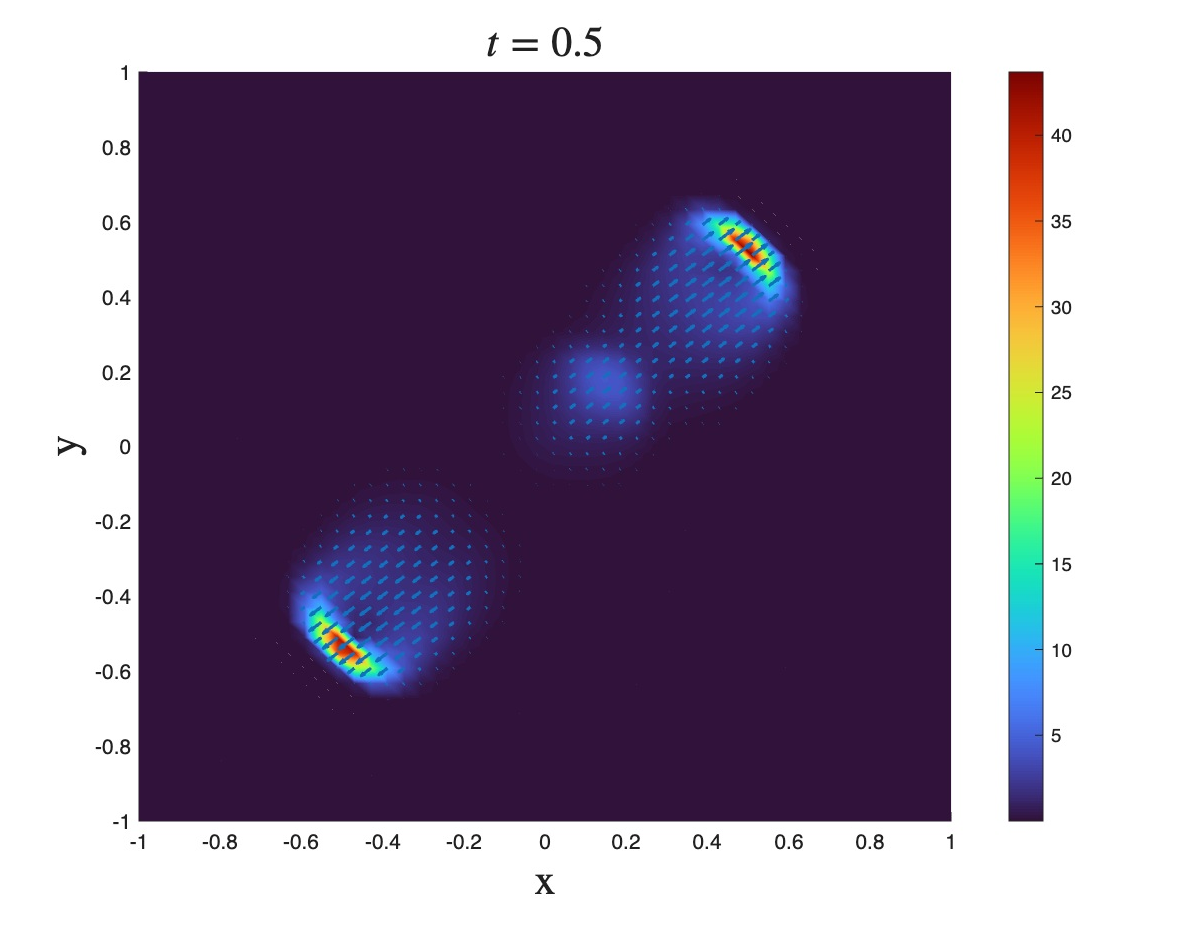}}
\subfigure{\includegraphics[scale=0.26]{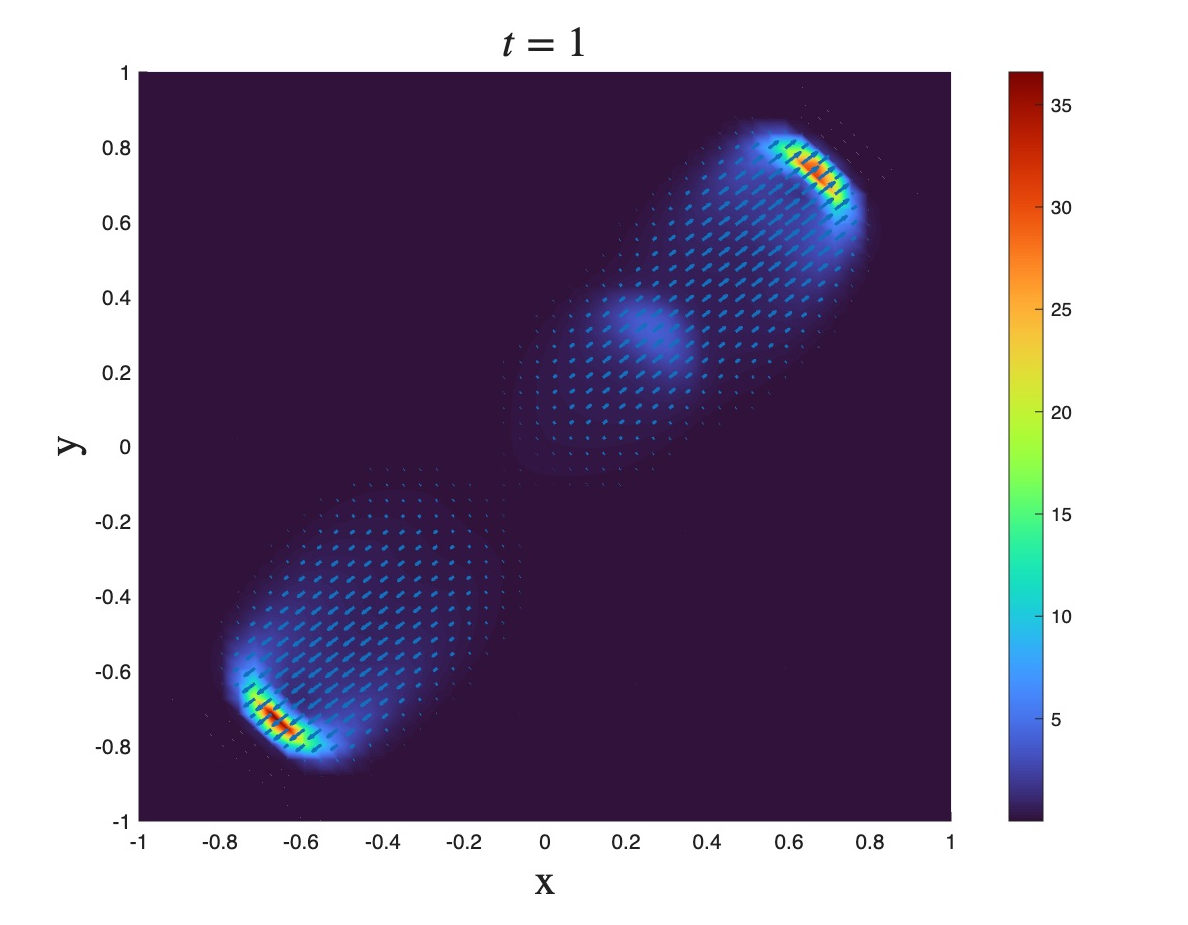}}
 \caption{Test 8: screenshots of the numerical simulation at (first line) the microscopic scale and (second line) the macroscopic scale.}
 
\label{Test_2D_positive}
\end{figure}

\paragraph{\textit{Tests 7-8: Sensitive dependence on initial data in $2$D}}
We then present a 2D scenario, generalizing the initial condition to $N$ agents, showing the sensitive dependence of the solutions on the initial data in higher dimension, both at microscopic and macroscopic scales.
First we run the microscopic model \eqref{sec2:CSmodel_pij} with $N=100$ agents. 
Test 7 and Test 8 differ on the initial condition.
At initial time,
two groups of 45 agents are randomly located in balls, centered in $(-0.3, -0.3)$ and $(0.3, 0.3)$, and radius $R=0.1$.
Each agent $i$ belonging to the group starts moving with initial velocity $V_i^0=(V_{ix}^0, V_{iy}^0)$ randomly chosen in $[-1, -0.2] \times [-1, -0.2] $, whereas initial velocities of agents belonging to the groups on the right are randomly chosen in $[0.2, 1] \times [0.2,1] $. 
 A small group of 10 agents with null velocities, is located in the ball centered in $(-0.08, -0.08)$ in Test 7, and in  $(0.08, 0.08)$ in Test 8. 
In the following $(X^0_i,Y^0_i)\in \mathbb{R}^{2}$  denotes the initial position of each agent $i$.
At the macroscopic scale, we run the system \eqref{eulereq} considering the spatial domain $[-1,1]\times[-1,1]$ in the time interval $[0,1]$, with a linear function $K$.
We approximate the solution implementing the scheme detailed in Section \ref{sec:schema_macro}, with $\Delta_x=\Delta_y=0.04$.
We consider an initial density given as the sum of Gaussian functions,
\begin{equation}
\rho_0(x,y)= \sum_{i=1}^{N}\frac{1}{2 \pi \sigma_x  \sigma_y} e^{-\frac{(x-X^0_i)^2}{2\sigma_x^2}-\frac{(y-Y^0_i)^2}{2\sigma_y^2}} 
\end{equation}
with $\sigma_x=\sigma_y=0.05$.
A plot of initial conditions for Test 7 and Test 8 is shown in Figure \ref{CI_2D} and Figure \ref{CI_2D_positive}, respectively.
Figure \ref{Test_2D} and Figure \ref{Test_2D_positive}  show several screenshots of the evolution of systems, both at microscopic and macroscopic scale. 
Depending on its initial position, the small group starts moving downward (Test 7) or upward (Test 8) reaching consensus with the closest group. A small variation on the initial data reflects in a different overall dynamics for the system, at both scales.
Note that the well-posedness issue does not depend on the number of agents (in the microscopic case, $N=3$ is anyway sufficient to observe this phenomenon). The evolution considered in Test 7 and Test 8 —with two clouds of particles moving in opposite directions—represents perfectly valid initial conditions, close to the pathological one, and serve as an example illustrating the sensitive dependence on the initial data.

\newpage

\section{Conclusion}\label{sec:Conclusion}
In this paper we investigated numerically flocking dynamics with topological interactions at different scales. The microscopic model represented the starting point of the analysis, for which several choices of interaction kernels were considered and analyzed. The corresponding Vlasov-type kinetic model was rigorously derived in \cite{BCR21} and its hydrodynamic limit to a pressureless Euler-type system was proved in \cite{BPR24}, under monokinetic assumptions on the initial condition. In our work we studied the agreement between different scales in several scenarios, also far from monokinetic initial conditions.
We finally focused on the sensitive dependence on the initial data in $1$D and $2$D scenarios, including the case of critical initial conditions leading to non uniqueness of solutions.
Future perspectives include further investigations of this aspect, including analogy with the selection of the limit for solutions to ODE with rough vector fields outside the Cauchy-Lipschitz theory. Moreover, a complete investigation of the macroscopic model, extending the analysis to the case of different topological interactions (e.g. attraction and repulsion effects) will be performed.  At kinetic and continuous scales, the role of topological interactions in continuous leader–follower transitions will be investigated, in the spirit of the work in \cite{CLMT_2025} where metric interactions are considered.
From a numerical point of view, we aim to extend the rescaling velocity method proposed in \cite{ReyTan} to the case of kinetic flocking models with topological interactions.

\section*{Acknowledgments}

 M.M. is funded by INdAM - GNCS Project, CUP E53C24001950001, entitled ``MODA:Integrating MOdel-based and DAta-Driven Methods for Multiscale Biological Systems''. T.T. is funded by the European Union’s Horizon Europe research and innovation programme under the Marie Skłodowska-Curie Doctoral Network DataHyking (Grant No.101072546).
 M.M. and T.T. are members of GNCS-INdAM research group.


\begin{thebibliography}{00}



\bibitem{ACH2018}
Giacomo Albi, Young-Pil Choi and Axel-Stefan Haeck. \textit{Pressureless Euler alignment system with control.} Mathematical Models and Methods in Applied Sciences 28.09 (2018): 1635-1664.

\bibitem{AF2024}
Giacomo Albi and Federica Ferrarese. \textit{Kinetic description of swarming dynamics with topological interaction and transient leaders}. Multiscale Modeling \& Simulation 22.3 (2024): 1169-1195. 

\bibitem{Ballerini}
Michele Ballerini, et al. \textit{Interaction ruling animal collective behavior depends on topological rather than metric distance: Evidence from a field study.} Proceedings of the national academy of sciences 105.4 (2008): 1232-1237.



\bibitem{Bialek}
William Bialek, et al. \textit{Statistical mechanics for natural flocks of birds.} Proceedings of the National Academy of Sciences 109.13 (2012): 4786-4791.

\bibitem{BCR21}
Dario Benedetto, Emanuele Caglioti and Stefano Rossi. \textit{Mean-field limit for particle systems with topological interactions.} Mathematics and Mechanics of Complex Systems 9.4 (2022): 423-440.

\bibitem{BD2016}
Adrien Blanchet and Pierre Degond. \textit{Topological interactions in a Boltzmann-type framework.} Journal of Statistical Physics 163 (2016): 41-60.

\bibitem{BPR24}
Dario Benedetto, Thierry Paul and Stefano Rossi. \textit{Propagation of chaos and hydrodynamic description for topological models.} Kinetic and Related Models 18.1 (2025): 19-34.

\bibitem{Carrillo2016}
José A. Carrillo, et al. \textit{Critical thresholds in 1D Euler equations with non-local forces.} Mathematical Models and Methods in Applied Sciences 26.01 (2016): 185-206.

\bibitem{Ciallella2023}
 Mirco Ciallella, Davide Torlo, and Mario Ricchiuto. \textit{Arbitrary high order WENO finite volume scheme with flux globalization for moving equilibria preservation.} Journal of Scientific Computing 96.2 (2023): 53.

\bibitem{CMPB21}
Emiliano Cristiani, Marta Menci, Marco Papi, Leonard Brafman. \textit{An all-leader agent-based model for turning and flocking birds.} Journal of Mathematical Biology 83.4 (2021): 45.

\bibitem{CLMT_2025}
Emiliano Cristiani, Nadia Loy, Marta Menci, Andrea Tosin. \textit{Kinetic description and macroscopic limit of swarming dynamics with continuous leader–follower transitions.} Mathematics and Computers in Simulation 228 (2025): 362-385.

 \bibitem{Haskovec2013} 
 Jan Haskovec. \textit{Flocking dynamics and mean-field limit in the Cucker–Smale-type model with topological interactions.} Physica D: Nonlinear Phenomena 261 (2013): 42-51.

\bibitem{cucker2007emergent}
Felipe Cucker and Steve Smale. \textit{Emergent behavior in flocks.} IEEE Transactions on automatic control 52.5 (2007): 852-862.

 \bibitem{cucker2007mathematics}
 Felipe Cucker and Steve Smale. \textit{On the mathematics of emergence.} Japanese Journal of Mathematics 2 (2007): 197-227.
 
 \bibitem{DP19}
 Pierre Degond and Mario Pulvirenti. \textit{Propagation of chaos for topological interactions.} The Annals of Applied Probability 29.4 (2019): 2594-2612.

 \bibitem{DPR23}
Pierre Degond, Mario Pulvirenti and Stefano Rossi. \textit{Propagation of chaos for topological interactions by a coupling technique.} Rendiconti Lincei 34.3 (2023): 641-655.

\bibitem{Gascon2001}
Llanos Gascón and José Miguel Corberán. \textit{Construction of second-order TVD schemes for nonhomogeneous hyperbolic conservation laws.} Journal of computational physics 172.1 (2001): 261-297.

\bibitem{KD}
Vijay Kumar and Rumi De. \textit{Efficient flocking: metric versus topological interactions.} Royal Society open science 8.9 (2021): 202158.

\bibitem{Kurganov2023}
 Alexander Kurganov, Yongle Liu, and Ruixiao Xin. \textit{Well-balanced path-conservative central-upwind schemes based on flux globalization.} Journal of Computational Physics 474 (2023): 111773.

\bibitem{KT}
Alexander Kurganov and Eitan Tadmor. \textit{Solution of two‐dimensional Riemann problems for gas dynamics without Riemann problem solvers.} Numerical Methods for Partial Differential Equations: An International Journal 18.5 (2002): 584-608.

\bibitem{LRS}
Daniel Lear, David N. Reynolds and Roman Shvydkoy. \textit{Global solutions to multi-dimensional topological Euler alignment systems.} Annals of PDE 8 (2022): 1-43.

\bibitem{LionsPaul}
Pierre-Louis Lions and Thierry Paul. \textit{Sur les mesures de Wigner.} Revista matemática iberoamericana 9.3 (1993): 553-618.

 \bibitem{MNP24}
 Marta Menci, Roberto Natalini and Thierry Paul. \textit{Microscopic, kinetic and hydrodynamic hybrid models of collective motions with chemotaxis: a numerical study.} Mathematics and Mechanics of Complex Systems 12.1 (2023): 47-83. 
 
 \bibitem{MPR25}
Marta Menci, Thierry Paul and Stefano Rossi. \textit{Propagation of chaos for topological models without regularity.} (2024). EMS Series of Congress Reports, In press.

\bibitem{TE_2024}
 Thierry Paul and Emmanuel Trélat. \textit{From microscopic to macroscopic scale equations: mean field, hydrodynamic and graph limits.} arXiv preprint arXiv:2209.08832 (2022).
 

\bibitem{ReyTan} 
Thomas Rey and Changhui Tan. \textit{An exact rescaling velocity method for some kinetic flocking models.} SIAM Journal on Numerical Analysis 54.2 (2016): 641-664.

\bibitem{RS}
 David N. Reynolds and Roman Shvydkoy. \textit{Local well-posedness of the topological Euler alignment models of collective behavior.} Nonlinearity 33.10 (2020): 5176.

\bibitem{ST_topo}
Roman Shvydkoy and Eitan Tadmor. \textit{Topologically based fractional diffusion and emergent dynamics with short-range interactions.} SIAM Journal on Mathematical Analysis 52.6 (2020): 5792-5839.

 \bibitem{Silverman}
Bernard W. Silverman, \textit{Density estimation for statistics and data analysis}, CRC press, 1986.

\bibitem{Sweby1984}
Peter K. Sweby. \textit{High resolution schemes using flux limiters for hyperbolic conservation laws.} SIAM Journal on Numerical Analysis 21, no. 5 (1984): 995-1011.

 \bibitem{Tan2017}
 Changhui Tan. \textit{A discontinuous Galerkin method on kinetic flocking models.} Mathematical Models and Methods in Applied Sciences 27.07 (2017): 1199-1221.

\bibitem{Terrell}
George R. Terrell and David W. Scott. \textit{Variable kernel density estimation.} The Annals of Statistics (1992): 1236-1265.
 
\bibitem{VZ}
 Tamás Vicsek and Anna Zafeiris. \textit{Collective motion.} Physics reports 517.3-4 (2012): 71-140.
 


 

\end{thebibliography}
\end{document}